\documentclass[10pt]{article}

\usepackage{latexsym,amsfonts}
\usepackage{amssymb,amsthm,upref,amscd}
\usepackage[T1]{fontenc}
\usepackage{times}

\newtheorem{Theorem}{Theorem}[section]
\newtheorem{definition}[Theorem]{Definition}
\newtheorem{lemma}[Theorem]{Lemma}
\newtheorem{proposition}[Theorem]{Proposition}

\newtheorem{remark}[Theorem]{Remark}
\newtheorem{Open Problem}[Theorem]{Open Problem}

\makeatletter \@addtoreset{equation}{section} \makeatother

\begin{document}

\title{\bf Existence  and concentration of  semiclassical states for nonlinear Schr\"odinger equations
 }

\author{Shaowei Chen   \thanks{E-mail adress:
chensw@amss.ac.cn (S. Chen),
 } \\ \\
\small  School of Mathematical Sciences, Capital Normal
University, \\
\small Beijing 100048, P. R. China\\
 }

\date{}
\maketitle

\begin{minipage}{13cm}
{\small {\bf Abstract:} In this paper, we study the following
semilinear Schr\"odinger equation
$$
-\epsilon^2\triangle u+ u+ V(x)u=f(u),\ u\in
H^{1}(\mathbb{R}^{N}),
$$
where $N\geq 2$ and $\epsilon>0$ is a small parameter. The
function $V$  is bounded in $\mathbb{R}^N,$
$\inf_{\mathbb{R}^N}(1+V(x))>0$ and it has a possibly degenerate
isolated critical point.  Under some conditions on  $f,$ we prove
that as $\epsilon\rightarrow 0,$ this equation has a solution
which concentrates at the
 critical point of $V$.

}
\end{minipage}

\renewcommand{\thefootnote}{}
\footnote{}

\footnote{}

\section{Introduction and main result}\label{8fhfgvvc1q23}
In this paper, we are concerned with  the following semilinear
Schr\"odinger equation
\begin{equation}\label{1.177qrreq}
-\epsilon^2\triangle u+ u+ V( x)u=f(u),\ u\in
H^{1}(\mathbb{R}^{N}),
\end{equation}
where $N\geq 2$  and $\epsilon>0$ is  a small parameter. The
function $f:\mathbb{R}\rightarrow \mathbb{R}$ satisfies
\begin{description}
\item{$(\bf{F_1}).$} $f\in C^1(\mathbb{R})$ and there exist $q\in
(2,2^*),$ $2<p_1<p_2<2^*$ and a constant $C>0$ such that
$$|f'(t)|\leq C(|t|^{p_1-2}+|t|^{p_2-2}),\ t\in\mathbb{R}$$
and for any $L>0$,
\begin{eqnarray}\label{l9d8d7f6}
\sup\{|f'(t)-f'(s)|/|t-s|^{q-2}\ |\ t,s\in[-L,L],\
t\neq s
 \}<\infty,\end{eqnarray} where   $2^*=2N/(N-2)$ if $N\geq 3$ and
$2^*=\infty$ if $N=2$;
\end{description}

\begin{description}
\item{$(\bf{F_2}).$} there exists $\mu>2$ such that $f(t)t\geq \mu
F(t)>0,$ $ t\neq 0$, where $F(t)=\int^t_0 f(s)ds$;
\end{description}

\begin{description}
\item{$(\bf{F_3}).$} $f(t)/|t|$ is an increasing function on
$\mathbb{R}\setminus\{0\}$;
\end{description}

\begin{remark}\label{hhat55ararr888audsd}
A typical function  which satisfies $\bf (F_1)-(F_3)$ is
$$f(t)=\sum^{m}_{i=1}a_i|t|^{\beta_i-2}t$$ with $2<\beta_1<\cdots<\beta_m<2^*$
and $a_i>0,$ $1\leq i\leq m.$
\end{remark}

The potential function  $V$ satisfies the following conditions:

\begin{description}
\item{$(\bf{V_0}).$} $\inf_{x\in\mathbb{R}^N}(1+V(x))> 0$ and
$\max_{x\in\mathbb{R}^N}|V(x)|<\infty;$
\end{description}

\begin{description}
\item{$(\bf{V_1}).$} $V\in C^{2}(\mathbb{R}^{N})$ has an isolated
critical point $x_0$ such that
$$V(x)=Q_{n^*}(x-x_0)+o(|x-x_0|^{n^*})$$ in some neighborhood of
$x_0,$ where $n^*\geq 2$ is an  even integer and $Q_{n^*}$ is an
$n^*$- homogeneous  polynomial in $\mathbb{R}^N$ which satisfies
that $\triangle Q_{n^*}\geq 0$ in $\mathbb{R}^N$ or $\triangle
Q_{n^*}\leq 0$ in $\mathbb{R}^N$ and $\triangle Q_{n^*}\not\equiv
0$ in $\mathbb{R}^N$.
\end{description}
\begin{remark}\label{jjfhvbb}
Without loss of generality, in what follows, we always assume that
$x_{0}=0.$ Typical examples for $Q_{n^*}$  are $\pm|x|^{n^*} $
$(n^*\geq 2).$
\end{remark}

Our main result of this paper is the following theorem
\begin{Theorem}\label{bbvv55erdff}
Suppose that $f$ satisfies $\bf (F_{1})-\bf (F_{4})$ and $V$
satisfies  $\bf (V_{0})$ and  $\bf (V_{1})$. Then there exist
$\epsilon_{0}>0$ and a set $\mathcal{K}$ whose elements are
radially symmetric solutions of equation
\begin{equation}\label{1.256de}
-\triangle u+ u=f(u),\ u\in H^{1}(\mathbb{R}^{N})
\end{equation}
 such that
 if $0<\epsilon<\epsilon_{0}$, then
  equation (\ref{1.177qrreq}) has a solution $u_{\epsilon}$   satisfying that
  $$\lim_{\epsilon\rightarrow 0}\mbox{dist}_{_{Y}}(v_{\epsilon},\mathcal{K})=0,$$
  where $v_{\epsilon}(x)=u_{\epsilon}(\epsilon x),$ $x\in\mathbb{R}^N$ and $Y=H^1(\mathbb{R}^N).$
\end{Theorem}

The analysis of the semilinear Schr\"odinger equation
(\ref{1.177qrreq}) has recently attracted a lot of attention due
to its many applications in mathematical physics.

If $v$ is a solution of equation (\ref{1.177qrreq}), then
$v(\epsilon x)$ is a solution of the following equation
\begin{equation}\label{1.1}
-\triangle u+ u+ V(\epsilon x)u=f(u),\ u\in H^{1}(\mathbb{R}^{N}).
\end{equation}   Equation (\ref{1.1}) is a
perturbation of the limit equation (\ref{1.256de}).
   If equation (\ref{1.256de})
  has  a solution $w\in
C^{2}(\mathbb{R}^{N})$   satisfying the  non-degeneracy condition:
$$\ker  L_0=\mbox{span}\left\{\frac{\partial \omega }{\partial x_i}\ |\ 1\leq
i\leq N\right\},$$ where $L_0 v=-\triangle v+v-f'(\omega)v$, then
in the celebrated paper \cite{Ambrosetti Badiale } (see also
\cite{Ambrosetti Secchi}), Ambrosetti, Badiale and Cingolani
developed a kind of variational reduction method and showed that
if the potential function $V$ has a strictly local minimizer or
maximizer $x_0$,
 then equation (\ref{1.1}) admits a solution
$u_\epsilon$ which converges to $\omega(\cdot-x_0)$ in
$H^{1}(\mathbb{R}^N)$ as $\epsilon\rightarrow 0$. In their
argument, the non-degeneracy property of $\omega$ plays essential
role.  Using the non-degeneracy condition and the reduction
method,  it was shown by Kang and Wei \cite{KW} that, at a strict
local maximum point $x_0$ of $V$ and for any positive integer $k$,
(\ref{1.177qrreq}) has a positive solution with $k$ interacting
bumps concentrating near $x_0$, while at a non-degenerate local
minimum point of $V(x)$ such solutions do not exist. Moreover,
under the assumption of the non-degeneracy condition, multiplicity
of solutions with one bump has also been considered by Grossi
\cite{Grossi}.

However, for a general nonlinearity $f$, it is very difficult to
verify the non-degeneracy condition for a solution of
(\ref{1.256de}). An effective method to attack problem
(\ref{1.177qrreq}) without using the non-degeneracy condition is
variational method. In \cite{Rabi2}, Rabinowitz used a global
variational method to show the existence of least energy solutions
for (\ref{1.177qrreq}) when $\epsilon>0$ is small, and the
condition imposed on $V$ is a global one, namely
$$0<\inf_{x\in\mathbb{R}^N}(1+V(x))<\liminf_{|x|\rightarrow\infty}(1+V(x)).$$
 In \cite{pino}, \cite{Del felmer},
\cite{Del felmer2}, \cite{Del felmer3} and \cite{Gui}, Del Pino,
Felmer and Gui used  different variational methods to obtain
nontrivial solution of (\ref{1.177qrreq}) for small $\epsilon>0$
under local conditions which can be roughly described as follows:
$V$ is local H\"older continuous on $\mathbb{R}^N,$
\begin{eqnarray}\label{mmfyftg}
\inf_{x\in \mathbb{R}^N}(1+V(x))> 0
\end{eqnarray}
and there exists $k$ disjoint bounded regions
$\Omega_1,\cdots,\Omega_k$ in $\mathbb{R}^N$ such that
\begin{eqnarray}\label{nnds66dtrf}
\inf_{x\in\partial\Omega_i}V(x)>\inf_{x\in\Omega_i}V(x).
\end{eqnarray}
 Their methods involve the deformation of
nonlinearity $f$ and some prior estimates. Recently, Byeon,
Jeanjean and Tanaka \cite{Beyon-Jeanjean} \cite{Beyon-Jeanjean2}
developed the variational methods and made great advance in
problem (\ref{1.177qrreq}).  Byeon and Jeanjean showed in
\cite{Beyon-Jeanjean} that if $N\geq3$, $V$ satisfies
(\ref{mmfyftg}) and (\ref{nnds66dtrf}) with $k=1$ and $f$
satisfies
\begin{description}\label{bbddrdf}
\item{$\bf (f_1).$} $f:\mathbb{R}\rightarrow\mathbb{R}$ is
continuous and $\lim_{t\rightarrow 0+}f(t)/t=0;$ \item{$\bf
(f_2).$} there exists some $p\in (1,2^*-1)$ such that
$\lim_{t\rightarrow\infty}f(t)/t^p<\infty;$\item{$\bf (f_3).$}
there exists $T>0$ such that $\frac{1}{2}mT^2<F(T)$, where
$F(t)=\int^{t}_{0}f(s)ds$ and $m=\inf_{x\in\Omega_1}V(x),$
\end{description}
then (\ref{1.177qrreq}) exists positive solution $v_\epsilon$
concentrating in the minimizers of $V$ in $\Omega_1$ as
$\epsilon\rightarrow 0.$ And in \cite{Beyon-Jeanjean2}, Byeon,
Jeanjean and Tanaka considered the case  $N=1,2$   and obtained
similar results. Their conditions on the nonlinearity $f$ are
almost optimal. Moreover, when $V$ satisfies (\ref{mmfyftg}) and
(\ref{nnds66dtrf}) with $k>1$ and $f$ satisfies $\bf (f_1)-(f_3)$,
in \cite{jeanjean3}, Cingolani, Jeanjean and Secchi constructed
multi-bump solutions for magnetic nonlinear Sch\"odinger equations
which contain equation (\ref{1.177qrreq}) as a special case.

  Comparing to the variational methods mentioned above,
   the Lyapunov reduction method of Ambrosetti and Badiale, although it
need the non-degeneracy condition, has its advantages that their
method can be used to deal with  elliptic equations involving
critical Sobolev exponent (see, for example, \cite{Ambrosetti
peral}) and other problems involving concentration compactness
(see, for example, \cite{Felli}).

In this paper, we indent to attack the problem (\ref{1.177qrreq})
though a Lyapunov reduction method, but avoiding  the
non-degeneracy condition for the solutions of limit equation
(\ref{1.256de}). In this paper, we develop a new  reduction method
for an isolated critical set $\mathcal{K}$ of the functional
corresponding to (\ref{1.256de}). This method can be regarded as a
generalization of Ambrosetti and  Badiale's method. The
non-degeneracy conditions for the solutions in this critical set
are no longer necessary and it does not involve the deformation of
nonlinearity. By combination of the new  reduction method and
Conley index theory which was developed by   Chang and Ghoussoub
in \cite{CG}(see also \cite{cha2}), we obtain a solution of
(\ref{1.1}) in a neighborhood of $\mathcal{K}$ for sufficiently
small $\epsilon>0.$ Our method is new and it can be used to other
problems which involve concentration compactness. In contrast with
the results of Byeon, Jeanjean and Tanaka, although the
assumptions we imposed on the nonlinearity $f$ are much stronger,
the assumptions we made on  $V$ seem weaker in a sense,    because
by the assumption $\bf (V_1)$, $x_0$ can be  a local maximum point
of $V.$

This paper is organized as follows: In section \ref{dcfrr5r}, we
obtain a critical  set of  the  functional  corresponding to
(\ref{1.256de}) with nontrivial Topology.  In section
\ref{ggdfsrseese11q} and section \ref{dadczd443eee}, a reduction
for the function corresponding to (\ref{1.1}) is developed. In
section \ref{ttrfd444ew3}, we give the proof of Theorem
\ref{bbvv55erdff}. Section \ref{gdg88futq9oi} and
\ref{99skdnfhfjfgevfd} are appendixes.

\bigskip

\noindent{\bf Notations.}  $\mathbb{R},\ \mathbb{Z}$ and
$\mathbb{N}$ denote the sets of real number, integer and positive
integer respectively. Let $E$ be a metric space. $B_{E}(a,\rho)$
denotes the open ball in $E$ centered at $a$ and having radius
$\rho$. The closure of a set $A\subset E$ is denoted by
$\overline{A}$ or $cl_E(A).$ $\mbox{dist}_E(a, A)$ denotes the
distance from  the point $a$ to the set $A\subset E$.  By
$\rightarrow$ we denote the strong and by $\rightharpoonup$ the
weak convergence.
  By $\ker A$
 denotes the null space of the operator $A.$ If $g$ is a $C^{2}$
 functional defined on a Hilbert
space $H$, $\nabla g$ (or $Dg$) and $\nabla^{2}g$ (or $D^2 g$)
denote the gradient of $g$ and the second derivative of $g$
respectively. And for $a, b\in\mathbb{R},$ we denote
$g^{a}:=\{u\in H\ | \ g(u)\leq a\}$ and $g_{b}:=\{u\in H\ | \
g(u)\geq b\}$ the sub- and super-level sets of the functional $g,
$ moreover, $g^{a}_{b}:=\{u\in H\ | \ b\leq g(u)\leq a\}.$
$\delta_{i,j}$ denotes the Kronecker notation, i.e.,
$\delta_{i,j}=1$ if $i=j$ and $0$ if $i\neq j.$ For a Banach space
$E,$ denote $\mathcal{L}(E)$ the Banach space consisting of all
bounded linear operator from $E$ to $E.$ If $H$ is a Hilbert space
and $W$ is a closed subspace of $H,$ we denote the orthogonal
complement space of $W$ in $H$ by $W^{\bot}.$ For a subset
$A\subset H,$ $\mbox{span}\{A\}$ denotes the subspace of $H$
generated by $A.$ For  a topology pair $(A,B)$ in metric space,
$\check{H}^*(A,B)$ denotes the
$\check{\mbox{C}}\mbox{ech}$-Alexander-Spanier cohomology with
coefficient group $\mathbb{Z}_2$ (see \cite{spanier}).

\section{Critical  sets of  limit functional with nontrivial Topology}\label{dcfrr5r}
 Throughout this paper, we denote the Sobolev
space $H^{1}(\mathbb{R}^{N})$ and  the radially symmetric function
space
$$H^{1}_{r}(\mathbb{R}^{N}):=\{u\in H^{1}(\mathbb{R}^{N})\ |\ u \
\mbox{is radially symmetric}\}$$ by  $Y$ and $X$ respectively.
 The inner product of
$Y$ is
$$\langle u,v\rangle=\int_{\mathbb{R}^{N}}(\nabla u\nabla v+
uv)dx,$$ and we use $||\cdot||$ to denote the norm of $Y$
corresponding  to this inner product. Define
$$I(u)=\frac{1}{2}\int_{\mathbb{R}^{N}}(|\nabla u|^{2}+ |u|^{2})dx
-\int_{\mathbb{R}^{N}}F(u)dx, \ u\in X.$$
$$J(u)=\frac{1}{2}\int_{\mathbb{R}^{N}}(|\nabla u|^{2}+ |u|^{2})dx
-\int_{\mathbb{R}^{N}}F(u)dx, \ u\in Y,$$
$$E_\epsilon(u)=\frac{1}{2}\int_{\mathbb{R}^{N}}(|\nabla u|^{2}
+ |u|^{2}+V(\epsilon x)|u|^{2})dx
-\int_{\mathbb{R}^{N}}F(u)dx, \ u\in Y.$$

For $h\in H^{-1}(\mathbb{R}^N)$, let $ (-\triangle+1)^{-1}h$ and
$(-\triangle+1+V(\epsilon x))^{-1}h$ be the solutions of
\begin{eqnarray}\label{uudgttdr77wew}
-\triangle u+u=h,\ u\in  H^{1}(\mathbb{R}^N)\end{eqnarray} and
\begin{eqnarray}\label{uudgttdrhhj77wew}
-\triangle u+u+V(\epsilon x)u=h,\ u\in  H^{1}(\mathbb{R}^N)
\end{eqnarray}
respectively.

Under conditions $\bf (F_1)-(F_3)$, $I$ satisfies Palais-Smale
condition (see, for example, \cite{Will}) and has a mountain pass
geometry, that is,
\begin{itemize}
\item[{\bf (i)}] $I(0)=0$, \item[{\bf (ii)}] there exist
$\rho_0>0$ and $\delta_0>0$ such that $I(u)\geq\delta_0$ for all
$||u||=\rho_0,$ \item[{\bf (iii)}] there exists $u_0\in X$ such
that $||u_0||>\rho_0$ and $I(u_0)<0.$
\end{itemize}
Thus the following minimax value is well defined and is larger
than $\delta_0,$
\begin{eqnarray}\label{xiquejiajja}
c=\inf_{\gamma\in\Gamma}\max_{t\in [0,1]} I(\gamma(t))
\end{eqnarray}
where
\begin{eqnarray}\label{66q5twfscccczz}
\Gamma=\{\gamma\in C([0,1],X)\ |\ \gamma(0)=0,\ I(\gamma(1))<0\}.
\end{eqnarray}

\begin{lemma}\label{yyrhft66t}For any $\sigma\in (0,\delta_0),$ if
$a\in (c-\sigma,c)$ and $b\in(c,c+\sigma)$ are regular values of
$I$, then $\check{H}^{1}(I^{b}, I^{a})\neq 0.$
\end{lemma}
\noindent{\bf Proof.} Since $b>c,$ by the definition of minimax
value $c,$ there exists $\gamma\in\Gamma$ such that
\begin{eqnarray}\label{hjdy66d5terffs555}
\max_{t\in[0,1]}I(\gamma(t))<b.
\end{eqnarray} Let $u_0=\gamma(1).$
We infer that $0$ and $u_0$ lie in different connected component
of $I^a.$ It follows that the homomorphism
$$\iota^*:\check{H}^0(I^a)\rightarrow \check{H}^0(\{0,u_0\})\cong \mathbb{Z}_2\oplus\mathbb{Z}_2$$
which is induced by the inclusion mapping
$\iota:\{0,u_0\}\hookrightarrow I^a$ is a surjection. Consider the
following homomorphism which is induced by the inclusion mapping
$j:\{0,u_0\}\hookrightarrow I^b,$
$$j^*:\check{H}^0(I^b)\rightarrow \check{H}^0(\{0,u_0\}).$$
By (\ref{hjdy66d5terffs555}), $0$ and $u_0$ lie in the same
connected component of $I^b.$  It follows that $j^*$ is not a
surjection.

Consider the following communicative diagram

\begin{center}
\setlength{\unitlength}{1cm}
\begin{picture}(7,3)

\put(-1,2){$\check{H}^{0}(I^{b})$}

\put(5.6,2){$\check{H}^{1}(I^{b},I^a)$}

\put(2.5,2){$\check{H}^{0}(I^a)$}

\put(0.4,2.1){\vector(1,0){1.8}}

\put(2.5,0){$\check{H}^{0}(\{0,u_0\})$}

\put(3,1.8){\vector(0,-2){1.3}}

\put(0.5,1.8){\vector(3,-2){2}}

 \put(1,2.3){\small $i^*$}

\put(3.3,1.05){\small $\iota^*$}

\put(0.7,1){\small $j^*$}

\put(4,2.1){\vector(1,0){1.5}}

\put(-2.5,2.1){\vector(1,0){1.4}}

\put(7.5,2.1){\vector(1,0){1.4}}

\put(4.7,2.3){\small $\alpha^*$}

\end{picture}
\end{center}
Since $j^*$ is not a surjection and $\iota^*$ is a surjection, by
this communicative diagram, we deduce that $\mbox{Image}(i^*)\neq
\check{H}^{0}(I^a).$  Moreover, by the property of exact sequence,
we have $\mbox{Image}(i^*)=\ker \alpha^*.$ Thus $\ker \alpha^*\neq
\check{H}^1(I^a)$. It follows that $\alpha^*\neq 0$.  Therefore,
$\check{H}^{1}(I^{b},I^a)\neq 0.$ \hfill$\Box$

\medskip

From Chapter 4 of \cite{Will}, we have the following lemma
\begin{lemma}\label{jjjag555sreesdbvv}
If  $ \nabla I(u)=0$ and $I(u)< 2c$, then $u$ does not change sign
in $\mathbb{R}^N$.
\end{lemma}

Let $\mathcal{F}$ be a $C^{1}$ functional defined on a Hilbert
space $M$
 with critical set
$K_{\mathcal{F}}$. And let $V$ be a pesudo-gradient vector field
 with respect to $D\mathcal{F}$ on $M$. A pesudo-gradient flow
associated with $V$ is the unique solution of the following
ordinary differential equation in $M:$
$$\dot{\eta}= -V(\eta(x,t)),\ \eta(x,0)=x.$$
 A subset $W$ of $M$ is
said to have the mean value property (for short (MVP)) if for any
$x\in M$ and any $t_{0}<t_{1}$ we have
$\eta(x,[t_{0},t_{1}])\subset W$ whenever  $\eta(x,t_{i})\in W,\
i=1,2.$

\begin{definition}\label{afsd556ggbdvf}(Definition I.10 of \cite{CG})
Let $\mathcal{F}$ be a $C^{1}$ functional on a Hilbert space $M$.
A subset $S$ of the critical set $K$ of $\mathcal{F}$ is said to
be a dynamically isolated critical set if there exist a closed
neighborhood $\mathcal{O}$ of $S$ and regular values $a<b$ of
$\mathcal{F}$ such that
\begin{equation}\label{homninchunenggqw545ff}
\mathcal{O}\subset \mathcal{F}^{-1}[a,b]
\end{equation}
and \begin{equation}\label{mnnavddf99kh}
cl(\widetilde{\mathcal{O}})\cap K\cap \mathcal{F}^{-1}[a,b]=S,
\end{equation}
where
$\widetilde{\mathcal{O}}=\bigcup_{t\in\mathbb{R}}\eta(\mathcal{O},t)$.
$(\mathcal{O}, a, b)$ is called an isolating triplet for $S.$
\end{definition}
\begin{definition}\label{haoyouduoooouagfdf}(Definition III.1
of \cite{CG})\label{zheqingsabainiaotefr66yh} Let $\mathcal{F}$ be
a $C^{1}$ functional on a Hilbet space $M$ and let $S$ be a subset
of the critical set $K_{\mathcal{F}}$ for $\mathcal{F}$. A pair
$(W,W_{-})$ of subset is said to be a GM pair for $S$ associated
with a pesudo-gradient vector field $V$, if the following
conditions hold:
\begin{description}
\item{(1).} $W$ is a closed (MVP) neighborhood of $S$ satisfying
$W\cap K=S$ and $W\cap \mathcal{F}_{\alpha}=\emptyset$ for some
$\alpha.$

\item{(2).} $W_{-}$ is an exit set for $W,$ i.e., for each
$x_{0}\in W$ and $t_{1}>0$ such that $\eta(x_{0},t_{1})\not\in W,$
there exists $t_{0}\in [0,t_{1})$ such that
$\eta(x_{0},[0,t_{0}])\subset W$ and $\eta(x_{0},t_{0})\in W_{-}.$

\item{(3).} $W_{-}$ is closed and is a union of a finite number of
sub-manifolds that transversal to the flow $\eta.$
\end{description}
\end{definition}

For $\alpha,\beta\in\mathbb{R},$ define
 $$\mathcal{K}^{\beta}_{\alpha}:=\{u\in X\ |\ \nabla
I(u)=0,\ \alpha\leq I(u)\leq \beta\}.$$ Let $a$ and $b$ are the
regular values which come from Lemma \ref{yyrhft66t}. Then by
Definition \ref{haoyouduoooouagfdf}, $\mathcal{K}^{b}_{a}$ is a
dynamically isolated critical set of $I$. By Lemma \ref{yyrhft66t}
and Theorem III.3
 of \cite{CG}, we have the following lemma
\begin{lemma}\label{yryhfgyhgtgftt645t} Let $\sigma>0$ be sufficiently small and $a\in(c-\sigma,c)$,
$b\in(c,c+\sigma)$ be regular values of $I$. If $(W,W_-)$ is a GM
pair of $\mathcal{K}^{b}_{a}$ associated with  some
pseudo-gradient vector field of $I$, then
$$\check{H}^{1}( W, W_-)\neq 0.$$
\end{lemma}
\begin{remark}\label{lfhjf998fu} In this remark, we shall show that
the set of regular values of $I$ is  dense in $\mathbb{R}.$
Therefore, for any  $\sigma>0$, there always exist regular values
of $I$ in $(c-\sigma,c)$ and $(c,c+\sigma).$ In fact, we shall
show that  $I(C)$ is of first category, where $C$ is the set of
critical points of $I$. It  suffices to prove that for any $u\in
C$, there exists $\delta_u>0$ such that $\overline{I(C\cap
B_X(u,\delta_u))}$ does not contain interior  points.

Let $u\in C$. Since $u$ is radially symmetric, the dimension of
the kernel space of the following operator is at most one
\begin{eqnarray}\label{i99o2w}
 \nabla^2 I(u):X\rightarrow X,\ h\in X\mapsto
h-(-\triangle+1)^{-1}f'(u)h.\nonumber
\end{eqnarray}

If $\dim \nabla^2 I(u)=0,$ then by Morse Lemma (see, e.g., Lemma
4.1 of \cite{cha}), there exists $\delta_u>0$ such that $u$ is the
unique critical point of $I$ in $B_X(u,\delta_u).$ Thus, in this
case, $I(C\cap B_X(u,\delta_u))=\{I(u)\}.$

If $\dim \nabla^2 I(u)=1,$ let $N=\mbox{ker}\nabla^2 I(u)$ and
note that $I$ is a $C^2$ functional, then by Lemma 1 of \cite{GM}
(see also Theorem 5.1 of \cite{cha}), there exist an origin
preserving $C^1$ diffeomorphism $\Phi$ of some $B_X(0,\delta_u)$
into $X$ and an an origin preserving $C^1$  map $h$ defined in
$N\cap B_X(0,\delta_u)$ into $X$ such that
$$I\circ\Phi(z,y)=I(u)+||Pz||^2-||(\mbox{id}-P)z||^2+I(h(y)+y)$$
where $P:N^\bot\rightarrow N^\bot$ is an orthogonal projection and
$N^\bot$ is the orthogonal complement of $N$ in $X$. Let $U=\{y\in
N\cap B_X(0,\delta_u)\ |\ h(y)+y\}.$ Then $U$ is a $C^1$
one-dimensional manifold.  Let us restrict $I$ to $U$. Then
$I:U\rightarrow\mathbb{R}$ is $C^1$. Moreover, $C\cap
B_X(0,\delta_u)=C\cap U,$ so $I(C\cap B_X(0,\delta_u))=I(C\cap
U)$. Therefore, by classical Sard theorem, $\overline{I(C\cap
B_X(0,\delta_u))}$ does not contain interior points.
\end{remark}

For  $r>0,$ $A\subset X,$ let
\begin{eqnarray}\label{bbfgr64qinniantt4trff}
N_{r}(A):=
 \{v\in X \ | \ \mbox{dist}_{X}(v,A)<r \}.\end{eqnarray}

 \begin{lemma}\label{bcng88yuyjhjnew55r}
 Let $c$ be the mountain pass value  coming from
Lemma \ref{yyrhft66t}. For any $r>0,$ there exists $\sigma_r>0$
such that if  $a\in(c-\sigma_r,c)$ and $b\in(c,c+\sigma_r)$ are
regular values of $I$, then
  there exists a  GM pair $(W,W_-)$ of
  the critical set
 $\mathcal{K}^{b}_{a}$ of the functional $I$ associated with
 the negative gradient vector field of $I$  such that $W\subset
N_{r}(\mathcal{K}^b_a).$
 \end{lemma}
\noindent{\bf Proof.} By $\bf (F_1)-(F_3),$ we know that $I$
satisfies   the Palais-Smale condition (see \cite{Will}).
Therefore, for any $r>0,$ there exists $\kappa_{r}>0$ such that if
$a\in(c-1,c)$ and $b\in(c,c+1)$, then
\begin{eqnarray}\label{gsn00osjll}
||\nabla I(v)||\geq \kappa_{r},\ \forall v\in I^{-1}[a,b]\setminus
N_{r/3}(\mathcal{K}^{b}_{a}).
\end{eqnarray}
Let
\begin{eqnarray}\label{bbctdgg0}
0<\sigma_r<\min\{r\kappa_{r}/6,1\}
\end{eqnarray}
and $a\in(c-\sigma_r,c)$ and $b\in(c,c+\sigma_r)$ be regular
values of $I$. For
\begin{eqnarray}\label{gdgb88fufhhhf776}
u\in I^{-1}[a,b] \cap N_{r/3}(\mathcal{K}^{b}_{a}),\end{eqnarray}
consider the negative gradient flow:
\begin{eqnarray}\label{bx888zz8z}
\dot{\eta}(t)=-\nabla I(\eta(t)),\ \eta(0)=u.\end{eqnarray}
Let
$$T^+_{u}=\sup\{t\geq 0\ |\ \mbox{for every}\ s\in[0,t],\
I(\eta(s))\geq a\}$$ and
$$T^-_{u}=\inf\{t\leq 0\ |\ \mbox{for every}\ s\in[t,0],\
I(\eta(s))\leq b\}.$$
 Let
$$U=\bigcup_{t\in [T^-_u,T^+_u]}\{\eta(t,u)
\ |\ u\in I^{-1}[a,b] \cap N_{r/3}(\mathcal{K}^{b}_{a})\}.$$ Then
$$[\mathcal{K}^{b}_{a}]\subset U,$$
where $$[\mathcal{K}^{b}_{a}] =\{v\in X\ | \
\omega(v)\cup\omega^*(v)\in \mathcal{K}^{b}_{a}\},$$
$\omega(v)=\cap_{t>0}\overline{\eta(v,[t,+\infty))}$ is the
$\omega-$limit set of $v$ and
$\omega^*(v)=\cap_{t>0}\overline{\eta(v,(-\infty,-t])}$ is the
$\omega^*-$limit set of $v$.

By \cite[Proposition III.2]{CG}, we deduce that there exists a
 GM pair $(W,W_-)$ of
$\mathcal{K}^{b}_{a}$ such that $W\subset U$. Thus, to prove this
Lemma, it suffices to prove that  if $\sigma_r>0$ is small enough,
then for  $u$ which satisfies (\ref{gdgb88fufhhhf776}),
\begin{eqnarray}\label{ncn99ysysh}
\sup_{t\in(T^-_{u},T^+_{u})}||\eta(t)-u||\leq \frac{2}{3}r.
\end{eqnarray}
Since their arguments are similar, we only give the proof for
\begin{eqnarray}\label{ncn996arrysysh}
\sup_{t\in[0,T^+_{u})}||\eta(t)-u||\leq \frac{2}{3}r.
\end{eqnarray}

 If (\ref{ncn996arrysysh}) were not true, then there
exist $0\leq t_1<t_2<T^+_{u}$ such that
$$r/3\leq||\eta(t)-u||\leq 2r/3,\ \forall t\in[t_1,t_2]$$
\begin{eqnarray}\label{h77dtdt}||\eta(t_1)-u||=r/3,\ ||\eta(t_2)-u||= 2r/3.
\end{eqnarray} According to (\ref{gsn00osjll}), we have
\begin{eqnarray}\label{vcbgftr655453g}
&&b-a\geq
I(\eta(t_1))-I(\eta(t_2))\nonumber\\
&&=\int^{t_1}_{t_2}\langle \nabla I(\eta(t)),\dot{\eta}(t)\rangle
dt=\int^{t_2}_{t_1}||\nabla I(\eta(t))||^2
dt\geq\kappa^{2}_{r}(t_2-t_1).\nonumber
\end{eqnarray}
It follows that
\begin{eqnarray}\label{bvn8846rttr}
t_2-t_1\leq (b-a)/\kappa^{2}_{r}.
\end{eqnarray}
Combining  (\ref{h77dtdt}) and (\ref{bvn8846rttr}) leads to
\begin{eqnarray}\label{bcn99udidii}
\frac{r}{3}&\leq&||\eta(t_2)-\eta(t_1)||\leq\int^{t_2}_{t_1}||\dot{\eta}(t)||dt\nonumber\\
&\leq&(t_2-t_1)^{1/2}(\int^{t_2}_{t_1}||\dot{\eta}(t)||^2)^{1/2}=
(t_2-t_1)^{1/2}(\int^{t_2}_{t_1}||\nabla I(\eta(t))||^2)^{1/2}\nonumber\\
&\leq&(t_2-t_1)^{1/2}(b-a)^{1/2}\leq
(b-a)/\kappa_{r}<2\sigma_r/\kappa_{r}.\nonumber
\end{eqnarray}
It contradicts (\ref{bbctdgg0}). Thus, (\ref{ncn996arrysysh})
holds. \hfill$\Box$

\section{A variational reduction for the limiting functional $I$}\label{ggdfsrseese11q}
 Let $\sigma>0$ be sufficiently small and
$a\in(c-\sigma,c)$, $b\in(c,c+\sigma)$ be regular values of $I$,
where  $c$ is defined by (\ref{xiquejiajja}). In what follows, for
the sake of  simplicity, we denote the critical set
$\mathcal{K}^{b}_{a}$ by $\mathcal{K}.$

 By \cite{Ber-Lions},  if
$u\in Y$ is a weak solution of
\begin{equation}\label{882rrrw} -\triangle u+ u
=f(u),
\end{equation} then $u$ and $\frac{\partial u}{\partial x_i},$ $1\leq i\leq N$ satisfy
 exponential decay at infinity.  As a consequence,   $\mathcal{K}$ is a compact
 subset of $W^{2,2}(\mathbb{R}^N).$ If $u\in
Y$ is a solution of equation (\ref{882rrrw}), then $\frac{\partial
u}{\partial x_{i}},$ $i=1,\cdots,N$  are the eigenfunctions for
the  eigenvalue problem
\begin{equation}\label{7yrttt}
-\triangle h+ h=f'(u)h.
\end{equation}

\begin{remark}\label{ggvcffdfdds} By
\cite[Theorem C. 3.4]{Simon}),  any eigenfunction  of the
eigenvalue problem (\ref{7yrttt}) satisfies exponential decay
 at infinity.
\end{remark}

 The argument in  \cite[Page 970-971]{Dancer} implies the following
 Lemma.
\begin{lemma}\label{trtgfgg}
Suppose that $u\in X$ is a solution of equation (\ref{882rrrw})
and it does not change sign in $\mathbb{R}^N$. If $v\in Y$ is a
solution of (\ref{7yrttt}) and satisfies
$$\left\langle v, \frac{\partial u}{\partial
x_{i}}\right\rangle=0,\ i=1,\cdots,N,$$ then $v\in X.$
\end{lemma}
\begin{remark}\label{qqqqi88abbvcvcv}
By Lemma \ref{jjjag555sreesdbvv}, we infer that if
$u\in\mathcal{K}$, then $u$ does not change sign in
$\mathbb{R}^N$.
\end{remark}

As it has been mentioned above,   $\mathcal{K}$ is a compact
subset in $W^{2,2}(\mathbb{R}^N)$. Thus for any $u\in\mathcal{K}$
and any $\varsigma>0$, there exists $\tau_u>0$ such that
\begin{eqnarray}\label{gdtdcfxscsb} \sum^{N}_{j=1}\left|\left|\frac{\partial
v}{\partial x_j}- \frac{\partial u}{\partial
x_j}\right|\right|<\varsigma,\ \forall v\in \mathcal{K}\cap
B_{X}(u,2\tau_u).
\end{eqnarray}
Therefore, we can choose a finite open sub-covering of
$\mathcal{K}$
\begin{eqnarray}\label{jjsgrrseds5srrrrs}
\mathcal{A}=\{B_{X}(u_i,\tau_{u_i})\ |\ i=1,\cdots,s\}
\end{eqnarray}
from the open covering $\{ B_{X}(u,\tau_{u})\ |\
u\in\mathcal{K}\}$. Let $\zeta\in C^{\infty}([0,+\infty))$ be such
that $0\leq\zeta(t)\leq 1$ for all $t,$ $\zeta(t)=1$ for
$t\in[0,1/2]$ and $\zeta(t)=0$ for $t\in[1,\infty).$ Let
$$\xi_i(u)=\frac{\zeta(||u-u_i||/\tau_{u_i})}
{\sum^{s}_{i=1}\zeta(||u-u_i||/\tau_{u_i})},\ 1\leq i\leq s.$$
Then $\{\xi_i\ |\ 1\leq i\leq s \}$ is a $C^\infty$ partition of
unity corresponding to the covering $\mathcal{A}$.

  For $u\in\mathcal{K},$ let
$$Y_u:=\{h\in X \ |\ \nabla^2I(u)h=0\},\ Z_u:=\mbox{span}\{\frac{\partial u}{\partial x_i}\ |\
1\leq i\leq N\}.$$
 Let\begin{equation}\label{gdhr7r6fgd}
\mathcal{Y}=\mbox{span}\{\cup^{s}_{i=1} Y_{u_i}\}.\end{equation}
Let \begin{equation}\label{zhonyaodehhy} q=\dim
\mathcal{Y}.\end{equation}

Let $\{e_1, e_2,\cdots,e_q\}$ be an orthogonal normal base of
 $\mathcal{Y}$. As mentioned in Remark \ref{ggvcffdfdds},   for every $1\leq n\leq q,$
$e_n\in W^{2,2}_{r}(\mathbb{R}^N)$ and $e_n$ satisfies exponential
decay at infinity.

  Let $\{e'_1,
  e'_2\cdots\}$ be an orthogonal normal base of
  $\mathcal{Y}^\bot,$ where $\mathcal{Y}^\bot$ is the orthogonal complement space
  of $\mathcal{Y}$ in $X$.
   From the appendix A of this paper,
for every $k\in\mathbb{N},$
   there exists
\begin{eqnarray}\label{oiidffdjhff4w3w3w}
E_k:= \{\tilde{e}_{j,k}\ |\ 1\leq j\leq k\},
\end{eqnarray}
 such that
\begin{itemize}
\item[{\bf (i)}] For every $k,$
   $E_k\subset X\cap W^{2,2}_{r}(\mathbb{R}^N)$ and $E_k\bot\mathcal{Y}$;
  \item[{\bf (ii)}]  Every $\tilde{e}_{j,k}$ satisfies  exponential decay  at
  infinity,
  $\langle\tilde{e}_{j,k},\tilde{e}_{j',k}\rangle=\delta_{j,j'}$
  and
  \begin{eqnarray}\label{gfhb7fyfteew}
   \sup_{1\leq j\leq k}||\tilde{e}_{j,k}-e'_j||\leq1/2^k.\nonumber
   \end{eqnarray}
  \end{itemize}
For every $k, $ denote
\begin{eqnarray}\label{oiidffdjhff4w3w3w}
X_k:= \mbox{span}\{E_k\}\oplus\mathcal{Y}.\nonumber
\end{eqnarray}
Let $P_{k}: X\rightarrow X_k$ and $P^\bot_k:X\rightarrow X^\bot_k$
be the orthogonal projections, where $X^\bot_k$ is the orthogonal
complement space of $X_k$ in $X.$ By the definition of $X_k$ and
the  properties $\bf (i)$ and $\bf (ii)$ mentioned above, we have
the following Lemma which is easy to prove.
\begin{lemma}\label{gdbcv66rtfe}
For every $h\in X,$
$\lim_{k\rightarrow\infty}||h-P_{k}h||=\lim_{k\rightarrow\infty}||P_{k}^\bot
h||=0.$
\end{lemma}

\begin{lemma}\label{bbfbupauuanznanwdse}
For any $r>0,$  there exists $l_r\in\mathbb{N}$  such that if
$k\geq l_r$, then for every $v\in N_r(\mathcal{K})$,
$P^\bot_k\nabla^2 I(v)|_{X^\bot_k}$ is invertible and
\begin{eqnarray}\label{gdij66rtg5w5qdaxz}
||(P^\bot_k\nabla^2
I(v)|_{X^\bot_k})^{-1}||_{\mathcal{L}(X^\bot_k)}\leq 2.\nonumber
\end{eqnarray}
\end{lemma}
\noindent{\bf Proof.} For $w\in X^{\bot}_{k}$,
\begin{eqnarray}\label{ssinati99ifufyfttdd}
P^\bot_k\nabla^2
I(v)w=w-P^\bot_k(-\triangle+1)^{-1}f'(v)w.\nonumber
\end{eqnarray}
Denote the operator $w\mapsto P^\bot_k(-\triangle+1)^{-1}f'(v)w$
by $A_{v,k}.$ If we can prove that
\begin{equation}\label{jj66w5zhonguokexue}
\limsup_{k\rightarrow\infty}\sup\{||A_{v,k}||_{\mathcal{L}(X^\bot_k)}\
|\ v\in N_r(\mathcal{K})\}=0,\end{equation} then the conclusion of
this Lemma follows. If (\ref{jj66w5zhonguokexue}) were not true,
we can choose $v_k\in N_r(\mathcal{K})$ and $w_k\in X_k^\bot$
with $||w_k||=1$, $k=1,2,\cdots,$ such that
\begin{eqnarray}\label{xifenliejjf7dy}
\limsup_{k\rightarrow\infty}||A_{v_k,k}w_k||>0.
\end{eqnarray}
Without loss of generality, we assume that $v_k\rightharpoonup
v_0$ in $X$ and $w_k\rightharpoonup w_0$ in $X$ as
$k\rightarrow\infty.$
 Since for any $2\leq p<2^*$, $X$
can be  compactly embedded  into the radially symmetric $L^p$
space
 (see, for example, \cite[Corollary 1.26]{Will})
$$L^{p}_{r}(\mathbb{R}^N):=\{u\in L^{p}(\mathbb{R}^N)\ |\ u\  \mbox{ is radially symmetric}\},$$
combining the condition $\bf (F_1)$, we can get that
$$\lim_{k\rightarrow\infty}\sup\{\int_{\mathbb{R}^N}|f'(v_k)w_k
h-f'(v_0)w_0 h |\ |\ h\in X,\ ||h||\leq 1\}=0.$$ It follows that
\begin{eqnarray}\label{bbfgf7fyyqinger}
\lim_{k\rightarrow\infty}||(-\triangle+1)^{-1}(f'(v_k)w_k-f'(v_0)w_0)||=0.
\end{eqnarray}
By (\ref{bbfgf7fyyqinger}) and Lemma \ref{gdbcv66rtfe}, we deduce
that $\lim_{k\rightarrow\infty}||A_{v_k,k}w_k||=0.$ But this
contradicts (\ref{xifenliejjf7dy}). \hfill$\Box$

\medskip

For $u\in \mathcal{K}$, denote $X_k\oplus Z_u$ by $W_{u,k} $ and
let $W_{u,k}^\bot$ be the orthogonal complement space of $W_{u,k}$
in $Y.$ Let $P_{W_{u_i,k}}:Y\rightarrow W_{u_i,k}$ and
$P_{W^{\bot}_{u_i,k}}:Y\rightarrow W^{\bot}_{u_i,k}$ be the
orthogonal projections.

\begin{lemma}\label{bbfgfyytrtrtrtrrwdse} Suppose that
$\kappa:=\max\{\tau_{u_i}\ |\ 1\leq i\leq s\}$ is sufficiently
small, where $\tau_{u_i}$ comes from (\ref{jjsgrrseds5srrrrs}).
Then there exist $C>0$ and $l_\kappa\in\mathbb{N}$ such that if
$k\geq l_\kappa$ and $v\in B_{X}(u_i,\tau_{u_i})$ for some $1\leq
i\leq s$, then $P_{W^{\bot}_{u_i,k}}\nabla^2
J(v)|_{W^{\bot}_{u_i,k}}$ is invertible and
\begin{eqnarray}\label{gdg5w5qdaxz}
||(P_{W^{\bot}_{u_i,k}}\nabla^2
J(v)|_{W^{\bot}_{u_i,k}})^{-1}||_{\mathcal{L}(W^{\bot}_{u_i,k})}\leq
C.
\end{eqnarray}

\end{lemma}
\noindent{\bf Proof.} We note that for $w\in W^{\bot}_{u_i,k}$,
\begin{eqnarray}\label{bqffdqnbv99ifufyfttdd}
P_{W^{\bot}_{u_i,k}}\nabla^2
J(v)w=w-P_{W^{\bot}_{u_{i},k}}(-\triangle+1)^{-1}f'(u)w.\nonumber
\end{eqnarray}
Since for any $p\in[2,2^*)$, $X$ can be compactly embedded  into
the radially symmetric $L^p$ space, by the condition $\bf (F_1)$,
we deduce that $w\mapsto
P_{W^{\bot}_{u_{i},k}}(-\triangle+1)^{-1}f'(v)w$ is a compact
operator. It follows that $P_{W^{\bot}_{u_i,k}}\nabla^2
J(v)|_{W^{\bot}_{u_i,k}}$ is a Fredholm operator with index zero.
Therefore, if we can prove that there exists $C>0$ which is
independent of $k$ such that, for sufficiently large $k,$
\begin{eqnarray}\label{gdgjjjacccofg5w5qdaxz}
||P_{W^{\bot}_{u_i,k}}\nabla^2
J(v)w||_{\mathcal{L}(W^{\bot}_{u_i,k})}\geq \frac{1}{C}||w||,\
\forall w\in W^{\bot}_{u_i,k},\ \forall v\in
B_X(u_i,\tau_{u_i})\nonumber
\end{eqnarray}
then the conclusion of this Lemma follows.

Without loss of generality, we assume that $u_i\equiv u_1$ and for
the sake of simplicity, we denote the operator
$P_{W^{\bot}_{u_1,k}}\nabla^2 J(v)|_{W^{\bot}_{u_1,k}}$ by
$H_{v,k}.$
 If  such $C>0$ does not
exist, then there exist sequences $\{\tau^k_{u_1}\}$,
$\{v_k\}\subset X$ and $\{w_k\}\subset Y$ such that
$\tau^k_{u_{1}}\rightarrow 0$ as $k\rightarrow\infty,$ $v_k\in
B_{X}(u_{1},\tau^k_{u_{1}})$, $w_k\in W^{\bot}_{u_{1},k}$,
$||w_k||=1$, $k=1,2,\cdots$ and
\begin{eqnarray}\label{gdb88dudjjdhab88au}
\lim_{k\rightarrow\infty}||H_{v_k,k}w_k||=0.\end{eqnarray}Passing
to a subsequence, we may assume that $w_k\rightharpoonup w_0$ in
$Y$ as $k\rightarrow\infty.$ By $\tau^k_{u_{1}}\rightarrow 0$ as
$k\rightarrow\infty$ and the assumption that $\{v_k\}\subset
B_{X}(u_{1},\tau^k_{u_{1}})$, we get that
\begin{eqnarray}\label{ffdvcy6dtdt}
\lim_{k\rightarrow\infty}||v_k-u_{1}||=0.
\end{eqnarray}
By  $w_k\in W^{\bot}_{u_{1},k}$ and $w_k\rightharpoonup w_0$ in
$Y$, we get that $w_0\bot X\oplus Z_{u_{1}}$. Combining the
condition $\bf(F_1)$, (\ref{ffdvcy6dtdt}) and the fact that
$w_k\rightharpoonup w_0$ in $Y$  leads to
\begin{eqnarray}\label{bbv8iisusgsfsf}
\lim_{k\rightarrow\infty}||(-\triangle+1)^{-1}(f'(v_k)w_k-f'(u_{1})w_k)||=0
\end{eqnarray}
and \begin{eqnarray}\label{bbv8iisusgsfsfiiju}
\lim_{k\rightarrow\infty}||(-\triangle+1)^{-1}(f'(u_1)w_k-f'(u_{1})w_0)||=0.
\end{eqnarray}
By (\ref{bbv8iisusgsfsfiiju}) and (\ref{bbv8iisusgsfsf}), we get
that
\begin{eqnarray}\label{bbvhhhvg8iisusgsfsf}
\lim_{k\rightarrow\infty}||(-\triangle+1)^{-1}(f'(v_k)w_k-f'(u_{1})w_0)||=0.
\end{eqnarray}
By Lemma \ref{gdbcv66rtfe}, we deduce that
\begin{eqnarray}\label{gdbcggdydttdtthubei}
\lim_{k\rightarrow\infty}||P_{W^{\bot}_{u_{1},k}}h- P_{(X\oplus
Z_{u_{_1}})^\bot}h||=0,\ \forall h\in Y,
\end{eqnarray}
where $P_{(X\oplus Z_{u_1})^\bot}:Y\rightarrow(X\oplus
Z_{u_{_1}})^\bot$ is the orthogonal projection.
 By (\ref{bbvhhhvg8iisusgsfsf}) and (\ref{gdbcggdydttdtthubei}), we get that
\begin{eqnarray}\label{nnfnfjf88f7fggg}
&&\lim_{k\rightarrow\infty}||P_{W^{\bot}_{u_{1},k}}((-\triangle+1)^{-1}f'(v_k)w_k)-P_{(X\oplus
Z_{u_{_1}})^\bot}((-\triangle+1)^{-1}f'(u_{1})w_0)||=0.
\end{eqnarray}
 By definition,
\begin{eqnarray}\label{bqqnbv99ifufyfttdd}
H_{v_k,k}w_k=w_k-P_{W^{\bot}_{u_{1},k}}(-\triangle+1)^{-1}f'(v_k)w_k.
\end{eqnarray}
By (\ref{nnfnfjf88f7fggg}) and the assumption
$\lim_{k\rightarrow\infty}||H_{v_k,k}w_k||=0$, we deduce that
$\{w_k\}$ is compact in $Y.$ Therefore,  $||w_k-w_0||\rightarrow0$
as $k\rightarrow\infty.$ It follows that $||w_0||=1,$ since
$||w_k||=1$ for every $k.$

Sending $k$ into infinity in the equality
(\ref{bqqnbv99ifufyfttdd}), by $w_0\in (X\oplus Z_{u_{_1}})^\bot,$
(\ref{gdb88dudjjdhab88au}) and (\ref{nnfnfjf88f7fggg}), we get
that
\begin{eqnarray}\label{n1171762ff}
P_{(X\oplus
Z_{u_{_1}})^\bot}(w_0-(-\triangle+1)^{-1}f'(u_{1})w_0)=0.
\end{eqnarray}
By $w_0\bot X$ and $u_{1}\in X,$ we have
\begin{eqnarray}\label{malkloa9au}
&&\langle
w_0-(-\triangle+1)^{-1}f'(u_{1})w_0,h\rangle\nonumber\\
&=& \langle
w_0,h\rangle-\langle(-\triangle+1)^{-1}f'(u_{1})h,w_0\rangle=0,\
\forall h\in X.
\end{eqnarray}
Since for any $h\in Z_{u_{1}},$ $$
 h-(-\triangle+1)^{-1}f'(u_{1})h=0,$$
 we get that
\begin{eqnarray}\label{miialkloa9au}
&&\langle
w_0-(-\triangle+1)^{-1}f'(u_{1})w_0,h\rangle\nonumber\\
&=&\langle h-(-\triangle+1)^{-1}f'(u_{1})h,w_0\rangle =0,\ \forall
h\in Z_{u_{1}}.
\end{eqnarray}
By (\ref{malkloa9au}) and (\ref{miialkloa9au}), we get that
\begin{eqnarray}\label{n11717fff62ff}
P_{X\oplus Z_{u_{_1}}}(w_0-(-\triangle+1)^{-1}f'(u_{1})w_0)=0.
\end{eqnarray}
By (\ref{n1171762ff}) and (\ref{n11717fff62ff}), we obtain
$$w_0-(-\triangle+1)^{-1}f'(u_{1})w_0=0,$$
that is, $w_0$ is an eigenfunction of (\ref{7yrttt}) with
$u=u_{1}\in \mathcal{K}.$ But $w_0$ satisfies $w_0\bot X\oplus
Z_{u_{1}}$ and $||w_0||=1.$ This contradicts Lemma \ref{trtgfgg}.
\hfill$\Box$

\medskip

 For
$v\in\cup^{s}_{i=1}B_{X}(u_i,\tau_{u_i}),$ let
\begin{eqnarray}\label{6erwtsfga55a}
\mathcal{T}_v=\mbox{span}\{\sum^{s}_{i=1}\xi_{i}(v)\frac{\partial
u_i}{\partial x_j}\ |\ 1\leq j\leq N\}.
\end{eqnarray}
  The
space $X_k\oplus \mathcal{T}_v$ is denoted by $E_{v,k}$. Let
$P_{E^{\bot}_{v,k}}:Y\rightarrow E^{\bot}_{v,k}$ be the orthogonal
projection.

\begin{lemma}\label{gsfs5s4srfsf}Suppose that
  $\kappa=\max\{\tau_{u_i}\ |\ 1\leq i\leq s\}$ is sufficiently
  small.
 Then there exist $C'>0$ and $ l_\kappa\in\mathbb{N}$
 such that if $k\geq l_\kappa$, then
 for every $v\in
\cup^{s}_{i=1}B_{X}(u_i,\tau_{u_i})$,  the operator
$P_{E^{\bot}_{v,k}}\nabla^2 J(v)|_{E^{\bot}_{v,k}}$ is invertible
and
\begin{eqnarray}\label{fdgtr6rtwfrkghh}
||(P_{E^{\bot}_{v,k}}\nabla^2
J(v)|_{E^{\bot}_{v,k}})^{-1}||_{\mathcal{L}(E^{\bot}_{v,k})}\leq
C'.\end{eqnarray}
\end{lemma}
\noindent{\bf Proof.} As the proof of Lemma
\ref{bbfgfyytrtrtrtrrwdse}, it suffices to prove that there exists
$C'>0$ which is independent of $k$ such that, for sufficiently
large $k,$
\begin{eqnarray}\label{jan99aiifdgtr6hrtwfrkghh}
||P_{E^{\bot}_{v,k}}\nabla^2
J(v)w||_{\mathcal{L}(E^{\bot}_{v,k})}\geq\frac{1}{C'}||w||,\
\forall w\in E^{\bot}_{v,k},\ \forall v\in
\cup^{s}_{i=1}B_{X}(u_i,\tau_{u_i}).\end{eqnarray}

 Without loss of
generality, we assume that $v\in B(u_1,\tau_{u_{_{1}}})$. Let
$P_{X_k}: Y\rightarrow X_k$ and $P_{\mathcal{T}_v}:Y\rightarrow
\mathcal{T}_v$ be orthogonal projections.
 For $h\in Y,$
\begin{eqnarray}\label{gdbc6s5srsss}
P_{E^{\bot}_{v,k}}h=h-P_{X_k}h-P_{\mathcal{T}_v}h,
\end{eqnarray}
and \begin{eqnarray}\label{ywtsfsfvxc5}
P_{\mathcal{T}_v}h=\sum^{N}_{j=1}\Big\langle
h,\sum^{s}_{i=1}\xi_i(v)\frac{\partial u_i}{\partial
x_j}\Big\rangle \frac{\sum^{s}_{i=1}\xi_i(v)\frac{\partial
u_i}{\partial x_j}}{||\sum^{s}_{i=1}\xi_i(v)\frac{\partial
u_i}{\partial x_j}||^2}.
\end{eqnarray}
Since $\{\xi_i\ |\ 1\leq i\leq s \}$ is a  partition of unity, we
get that for every $1\leq j\leq N,$
\begin{eqnarray}\label{sfs55s4wf}
||\frac{\partial u_1}{\partial
x_j}-\sum^{s}_{i=1}\xi_i(v)\frac{\partial u_i}{\partial x_j}||
&=&||\sum^{s}_{i=1}\xi_i(v)\frac{\partial u_1}{\partial
x_j}-\sum^{s}_{i=1}\xi_i(v)\frac{\partial u_i}{\partial
x_j}||\nonumber\\
&\leq&\sum^{s}_{i=1}\xi_i(v)||\frac{\partial u_1}{\partial
x_j}-\frac{\partial u_i}{\partial x_j}||.
\end{eqnarray}
If $\xi_i(v)\neq 0,$ then $v\in B_{X}(u_i,\tau_{u_{i}})$.
Combining
 the assumption $v\in B_{X}(u_1,\tau_{u_{1}})$, we get that
$u_1\in B_{X}(u_i,2\tau_{u_{i}})\cap \mathcal{K}$. Therefore, by
(\ref{gdtdcfxscsb}), we deduce that
\begin{eqnarray}\label{bcbudydhggdf7yadui}
\sum^{s}_{i=1}||\frac{\partial u_1}{\partial x_j}-\frac{\partial
u_i}{\partial x_j}||<\varsigma,\ \mbox{if}\ \xi_i(v)\neq 0.
\end{eqnarray}
Combining (\ref{sfs55s4wf}) and (\ref{bcbudydhggdf7yadui}) leads
to
\begin{eqnarray}\label{nnfbch66dyquashidin}
||\frac{\partial u_1}{\partial
x_j}-\sum^{s}_{i=1}\xi_i(v)\frac{\partial u_i}{\partial
x_j}||<\varsigma,\ \mbox{for every }\ 1\leq j\leq N.
\end{eqnarray}
 Thus, there exists $C>0$ which is independent of $k$ such that
\begin{eqnarray}\label{4gdvc66dtdfaa}
||P_{\mathcal{T}_v}h-P_{Z_{u_1}}h||\leq C\varsigma ||h||,\ \forall
h\in Y,
\end{eqnarray}
where $$P_{Z_{u_1}}:Y\rightarrow Z_{u_1},\ h\mapsto
\sum^{N}_{j=1}\Big\langle h,\frac{\partial u_1}{\partial
x_j}\Big\rangle \frac{\frac{\partial u_1}{\partial
x_j}}{||\frac{\partial u_1}{\partial x_j}||^2} $$ is orthogonal
projection. By (\ref{gdbc6s5srsss}) and (\ref{4gdvc66dtdfaa}), we
have
\begin{eqnarray}\label{gdvc66dtdfaa}
||P_{E^{\bot}_{v,k}}h-P_{W^{\bot}_{u_{_{1}},k}}h||\leq C\varsigma
||h||,\ \forall h\in Y.
\end{eqnarray}
For $w\in E^{\bot}_{v,k},$ we have
\begin{eqnarray}\label{fsr5srsfalll}
&&||P_{E^{\bot}_{v,k}}\nabla^2 J(v)w||\\
&\geq& ||P_{W^{\bot}_{u_1,k}}\nabla^2
J(v)w||-||(P_{E^{\bot}_{v,k}}-P_{W^{\bot}_{u_1,k}})\nabla^2
J(v)w||\nonumber\\
&\geq&||P_{W^{\bot}_{u_1,k}}\nabla^2 J(v)w||-C\varsigma||\nabla^2
J(v)||_{\mathcal{L}(Y)}||w||\ (\mbox{by}\ (\ref{gdvc66dtdfaa}))\nonumber\\
&\geq&||P_{W^{\bot}_{u_1,k}}\nabla^2 J(v)(w-P_{Z_{u_1}}w)||-||
P_{W^{\bot}_{u_1,k}}\nabla^2
J(v)(P_{Z_{u_1}}w)||-C\varsigma||\nabla^2
J(v)||_{\mathcal{L}(Y)}||w||\nonumber\\
&\geq&C||w-P_{Z_{u_1}}w||-||\nabla^2
J(v)||_{\mathcal{L}(Y)}||P_{Z_{u_1}}w||\nonumber\\
&&-C\varsigma||\nabla^2 J(v)||_{\mathcal{L}(Y)}||w||\ (\mbox{by}\
w-P_{Z_{u_1}}w\in W^{\bot}_{u_1,k}\ \mbox{and}\
(\ref{gdg5w5qdaxz}))\nonumber\\
&\geq&C||w||-(C+||\nabla^2
J(v)||_{\mathcal{L}(Y)})||P_{Z_{u_1}}w||-C\varsigma||\nabla^2 J(v)||_{\mathcal{L}(Y)}||w||\nonumber\\
&=&C||w||-(C+||\nabla^2
J(v)||_{\mathcal{L}(Y)})||P_{\mathcal{T}_v}w-P_{Z_{u_1}}w||\nonumber\\
&&-C\varsigma||\nabla^2 J(v)||_{\mathcal{L}(Y)}||w||\ (\mbox{since}\ P_{\mathcal{T}_v}w=0)\nonumber\\
&\geq& C||w||-\varsigma C(C+||\nabla^2
J(v)||_{\mathcal{L}(Y)})||w||-C\varsigma||\nabla^2
J(v)||_{\mathcal{L}(Y)}||w||.\ (\mbox{by}\
(\ref{4gdvc66dtdfaa}))\nonumber
\end{eqnarray}
It follows that if $\kappa>0$ is  sufficiently small, then there
exist $l_\kappa\in\mathbb{N}$ and $C'>0$ such that for every
$k\geq l_\kappa,$ (\ref{jan99aiifdgtr6hrtwfrkghh}) holds.
\hfill$\Box$

\medskip

Recall that   $X^\bot_k$ is the orthogonal complement space of
$X_k$ in $X$ and  $P_{k}: X\rightarrow X_k$,
$P^\bot_k:X\rightarrow X^\bot_k$ are  orthogonal projections. Let
$$\mathcal{N}_{\delta,\tau,k}:=\{u+v\in X\
|\ u\in X_k,\ \mbox{dist}_X(u,P_{k}\mathcal{K})<\delta,\ v\in
X^{\bot}_{k} ,\  ||v||<\tau\},$$ where $P_{k}\mathcal{K}=\{P_k v\
|\ v\in\mathcal{K}\}.$ By Lemma \ref{gdbcv66rtfe} and the fact
that $\mathcal{K}$ is a compact subset of $X,$ we get that as
$k\rightarrow\infty$, the Hausdorff distance of $\mathcal{K}$ and
$P_k \mathcal{K}$,
\begin{eqnarray}\label{shengdanweirenjie}\sup_{v\in P_k
\mathcal{K}}\mbox{dist}_X(v, \mathcal{K}) +\sup_{u\in
\mathcal{K}}\mbox{dist}_X(u, P_k\mathcal{K})\rightarrow
0.\end{eqnarray} Thus, for any $\delta>0$, $\tau>0$ and
$0<r<\min\{\delta,\tau\}$, if $k$ is sufficiently large, then
\begin{equation}\label{gdteraa5wdd}
N_{r}(\mathcal{K})\subset \mathcal{N}_{\delta,\tau,k},
\end{equation}
where $N_{r}(\mathcal{K})$ comes from
(\ref{bbfgr64qinniantt4trff}). And for any $r>0$, if
$\delta,\tau\in(0,r/2)$, then for sufficiently large $k$,
\begin{equation}\label{gdtaazeraa5wdd}
 \mathcal{N}_{\delta,\tau,k}\subset N_{r}(\mathcal{K}).
\end{equation}

Let
\begin{eqnarray}\label{nnnzbxccc7y}
\mathcal{N}_{\delta,k}:=\{u\in X_k\ |\ \mbox{dist}_X(u,P_k
\mathcal{K})<\delta\}.
\end{eqnarray}

\begin{lemma}\label{77eygdjjj876yyyyff5ttt}
If    $\delta>0$ is sufficient small and $k$ is sufficiently
large, then there exists a $C^{1}-$mapping
$$\pi_{k}:\mathcal{N}_{\delta,k} \rightarrow  X^\bot_k, $$
satisfying
\begin{itemize}
\item[{\bf (i)}] $\langle \nabla I(v+\pi_{k}(v)), \phi\rangle=0,$
$\forall \phi\in X^{\bot}_k;$

 \item[{\bf (ii)}] $\lim_{k\rightarrow\infty}\sup\{||\pi_k(v)|| \ |\ v\in \mathcal{N}_{\delta,k}\}=0;$

 \item[{\bf (iii)}] $\lim_{k\rightarrow\infty}\sup\{||D\pi_k(v)h||
 \ |\ v\in \mathcal{N}_{\delta,k},
 \ h\in X_k,\ ||h||=1\}=0;$

  \item[{\bf (iv)}]  If $v$ is a critical point of $I(v+\pi_k(v))$,
  then $v+\pi_k(v)$ is a critical point of $I.$
\end{itemize}
\end{lemma}
\noindent{\bf Proof.} By Lemma \ref{bbfbupauuanznanwdse}, if $r>0$ is small enough, then the operator
$$L_{v,k}:=P^\bot_k\nabla^2 I(v)|_{X^\bot_k}: X^{\bot}_{k}\rightarrow X^{\bot}_{k}$$
is invertible and if $k\geq l_\kappa$,
\begin{eqnarray}\label{nnb99ifufjjj}
||L_{v,k}^{-1}||_{\mathcal{L}(X_k^\bot)}\leq 2, \ \forall v\in
N_{r}(\mathcal{K}).
\end{eqnarray}

Assume that $0<\delta<r$,  by (\ref{gdtaazeraa5wdd}), if $k$ is
large enough, then $\mathcal{N}_{\delta,k}\subset
N_{r}(\mathcal{K})$.

For $\rho>0$ and $v\in \mathcal{N}_{\delta,k},$ define
$$\Psi_{v,k}:\overline{B_{X^{\bot}_{k}}(0,\rho)}\rightarrow X^{\bot}_{k},\
w\mapsto w-L^{-1}_{v,k}P^{\bot}_{k}\nabla I(v+w).$$ For any
$w_i\in \overline{B_{X^{\bot}_{k}}(0,\rho)},$ $i=1,2,$ by  the
definition of  $L_{v,k}$, we have
$w_2-w_1-L^{-1}_{v,k}P^{\bot}_{k}\nabla^2 I(v)(w_2-w_1)=0.$
Therefore,
\begin{eqnarray}\label{ggdbchhhfjfj99ij}
&&||\Psi_{v,k}(w_2)-\Psi_{v,k}(w_1)||\nonumber\\
&=&||w_2-w_1-L^{-1}_{v,k}P^{\bot}_{k}\nabla^2 I(v+\theta
w_2+(1-\theta)w_1)(w_2-w_1)||\nonumber\\
&&(\mbox{by the mean value theorem},\ 0<\theta=\theta(x)< 1) \nonumber\\
&\leq& ||w_2-w_1-L^{-1}_{v,k}P^{\bot}_{k}\nabla^2
I(v)(w_2-w_1)||\nonumber\\
&&+||L^{-1}_{v,k}P^{\bot}_{k}(\nabla^2 I(v+\theta
w_2+(1-\theta)w_1)-\nabla^2 I(v))(w_2-w_1)||\nonumber\\
&=&||L^{-1}_{v,k}P^{\bot}_{k}(\nabla^2 I(v+\theta
w_2+(1-\theta)w_1)-\nabla^2 I(v))(w_2-w_1)||\nonumber\\
&\leq& 2||(\nabla^2 I(v+\theta w_2+(1-\theta)w_1)-\nabla^2
I(v))(w_2-w_1)|| \ (\mbox{by}\ (\ref{nnb99ifufjjj})).
\end{eqnarray}
Since $I\in C^2(X,\mathbb{R})$ and   $\mathcal{K}$ is compact in
$X$, if $\delta$ and $\rho$ are small enough, then for any
$v\in\mathcal{N}_{\delta,k}$ and $w\in
\overline{B_{X^{\bot}_{k}}(0,\rho)},$
$$||\nabla^2 I(v+w)-\nabla^2
I(v)||_{\mathcal{L}(X)}<1/4.$$ Thus,  by (\ref{ggdbchhhfjfj99ij}),
we get that for any $w_i\in \overline{B_{X^{\bot}_{k}}(0,\rho)},$
$i=1,2,$
\begin{eqnarray}\label{hhdhfgtdtdfdf}
||\Psi_{v,k}(w_2)-\Psi_{v,k}(w_1)||\leq\frac{1}{2}||w_2-w_1||.
\end{eqnarray}
If $\delta>0$ is  small enough and  $k$ is large enough, then for
every $v\in\mathcal{N}_{\delta,k}$,
$$||\Psi_{v,k}(0)||\leq\rho/2.$$
Then by (\ref{hhdhfgtdtdfdf}), we get that for every
$w\in\overline{B_{X^{\bot}_{k}}(0,\rho)}$,
\begin{eqnarray}\label{mmvnviiidy}
||\Psi_{v,k}(w)||\leq||\Psi_{v,k}(w)-\Psi_{v,k}(0)||+||\Psi_{v,k}(0)||\leq\rho.
\end{eqnarray}
By (\ref{hhdhfgtdtdfdf}) and (\ref{mmvnviiidy}),  $\Psi_{v,k}$ is
a contractive mapping in $\overline{B_{X^{\bot}_{k}}(0,\rho)}$ if
$\delta$ and $\rho$ are small enough and  $k$ is large enough.
Thus, by Banach fixed point theorem, there exists unique fixed
point $\pi_k(v)\in \overline{B_{X^{\bot}_{k}}(0,\rho)}.$ It is
easy to verify that $\pi_k$ is a $C^1-$mapping and it satisfies
the result $\bf (i)$.

Now, we  give the proof of $\bf (ii).$ By $P_k^\bot\nabla
I(v+\pi_k(v))=0$ and $\pi_k(v)\in X^{\bot}_{k}$, we get that
\begin{eqnarray}\label{bcv77yeter}
0&=&\langle\nabla I(v+\pi_k(v)),\pi_k(v)\rangle\nonumber\\
&=&||\pi_k(v)||^2-\int_{\mathbb{R}^N}f(v+\pi_k(v))\cdot\pi_k(v).
\end{eqnarray}
By Lemma \ref{gdbcv66rtfe}, we deduce that for any sequence
$\{v_k\}$ with $v_k\in\mathcal{N}_{\delta,k}$,
$\pi_k(v_k)\rightharpoonup 0$ in $X$ as $k\rightarrow\infty$.
Combining  the compact embedding $X\hookrightarrow
L^{p}_{r}(\mathbb{R}^N)$, we obtain
$$\lim_{k\rightarrow\infty}\int_{\mathbb{R}^N}|f(v_k+\pi_k(v_k))|\cdot|\pi_k(v_k)|
=0.$$ It follows that
\begin{eqnarray}\label{nnvmvoofjfhfttg}
\lim_{k\rightarrow\infty}\sup\{\int_{\mathbb{R}^N}f(v+\pi_k(v))\cdot\pi_k(v)\
| \ v\in \mathcal{N}_{\delta,k}\}=0.
\end{eqnarray}
 The conclusion $\bf (ii)$ follows from  (\ref{bcv77yeter}) and
 (\ref{nnvmvoofjfhfttg}).

Differentiating  equation $P^{\bot}_{k}\nabla I(v+\pi_k(v))=0$ for
the variable $v$ in the direction $h\in X_k$, we get that
 \begin{eqnarray}\label{hhdgdgdfdf77y}
D\pi_k(v)
 h-P^\bot_k(-\triangle+1)^{-1}f'(v+\pi_k(v))(h+D\pi_k(v)h)=0.
 \end{eqnarray}
Note that $D\pi_k(v)h\in X^{\bot}_{k}$. By (\ref{nnb99ifufjjj}),
(\ref{hhdgdgdfdf77y}) and $\lim_{k\rightarrow\infty}||\pi_k(v)||=
0$, we get that if $k$ is large enough, then
\begin{eqnarray}\label{mmdbdgdgdfdfff}
\frac{1}{2}||D\pi_k(v)h||&\leq& ||D\pi_k(v)
 h-P^\bot_k(-\triangle+1)^{-1}f'(v+\pi_k(v))D\pi_k(v)h||\\
 &=&||P^\bot_k(-\triangle+1)^{-1}f'(v+\pi_k(v))h||\nonumber
\end{eqnarray}
It follows that for sufficiently large $k,$
\begin{eqnarray}\label{ggdtdref55001223}
\sup\{||D\pi_k(v)h||\ |\ v\in\mathcal{N}_{\delta,k},\ h\in X_k,\
||h||\leq 1\}<\infty.\end{eqnarray} By (\ref{hhdgdgdfdf77y}), we
get that
\begin{eqnarray}\label{99fhfgtte5erre}
||D\pi_k(v)h||^2=\int_{\mathbb{R}^N}f'(v+\pi_k(v))\cdot(h+D\pi_k(v)h)\cdot
D\pi_k(v)h.\end{eqnarray}
 (\ref{ggdtdref55001223}) and
the same argument  as  (\ref{nnvmvoofjfhfttg}) yield
$$\lim_{k\rightarrow\infty}\sup\{\int_{\mathbb{R}^N}f'(v+\pi_k(v))\cdot(h+D\pi_k(v)h)\cdot
D\pi_k(v)h\ |\ v\in \mathcal{N}_{\delta,k},\ h\in X_k,\ ||h||\leq
1\}=0.$$ Combining  (\ref{99fhfgtte5erre}), we get the conclusion
$\bf (iii).$

By $\bf (iii),$ if $k$ is sufficiently large, then
$$\{h+D\pi_k(v)h\ |\ h\in X_k \}+ X^{\bot}_{k}=X.$$ Combining the result $\bf (i)$, we get that
 if $v_0$ is a
critical point of $I(v+\pi_k(v))$, then $v_0+\pi_k(v_0)$ is a
critical point of $I$. \hfill$\Box$

\begin{remark}\label{vvxfdrrd11887}
By $\bf (ii)$ and $\bf (iv)$ of Lemma
\ref{77eygdjjj876yyyyff5ttt}, $\mathcal{N}_{\delta,\tau,k}$ is a
neighborhood of $\mathcal{K}$ if \begin{eqnarray}\label{fsrsfsfsf}
\tau>\sup\{||\pi_k (v)|| \ |\ v \in \mathcal{N}_{\delta,k}
\}.\end{eqnarray}
\end{remark}

\begin{lemma}\label{yryuqsaxxx3et}  Let $\mathcal{I}_k(u)=\frac{1}{2}||P^\bot_k u||^2+I(P_k
u+\pi_{k}(P_k u )).$ Then $$
\lim_{k\rightarrow\infty}||\mathcal{I}_k-I||_{C^1(\overline{\mathcal{N}_{\delta,\tau,k}})}=0.$$
\end{lemma}
\noindent{\bf Proof.} By definition, we have
$$\mathcal{I}_k(u)=\frac{1}{2}|| u||^2+\frac{1}{2}||\pi_{k}(P_k u )||^2
-\int_{\mathbb{R}^N}F(P_k u+\pi_{k}(P_k u )).$$ For any sequence
$\{u_k\}$ with $u_k\in \overline{\mathcal{N}_{\delta,\tau,k}}$, by
the mean value theorem, we get that
\begin{eqnarray}\label{gdbvbyf6rtrgfd}
F(P_k u_k+\pi_{k}(P_k u_k ))-F(u_k)& =&\zeta(u_k,\theta)(P_k
u_k+\pi_{k}(P_k u_k )-u_k)\nonumber\\
&=&\zeta(u_k,\theta)(\pi_{k}(P_k u_k )-P_k^\bot u_k)\nonumber
\end{eqnarray}
where
\begin{eqnarray}\label{gdgbc77fyfgdddd}
\zeta(u_k,\theta)= f'(\theta P_k u_k+\theta\pi_{k}(P_k u_k
)+(1-\theta)u_k)\nonumber
\end{eqnarray}
with $0<\theta(x)<1,$ $x\in \mathbb{R}^N$. Then we have
\begin{eqnarray}\label{hfbv8fyfyyfww}
\int_{\mathbb{R}^N}\Big|F(P_k u_k+\pi_{k}(P_k u_k
))-F(u_k)\Big|=\int_{\mathbb{R}^N}|\zeta(u_k,\theta)|\cdot|\pi_{k}(P_k
u_k )-P_k^\bot u_k|.
\end{eqnarray}
By $\bf (ii)$ of Lemma \ref{77eygdjjj876yyyyff5ttt}, we get that
for every $2\leq p<2^*,$
\begin{eqnarray}\label{bcbvvvvv55}
\lim_{k\rightarrow\infty}\int_{\mathbb{R}^N}|\pi_{k}(P_k u_k
)|^{p}=0.
\end{eqnarray}
By Lemma \ref{gdbcv66rtfe},we have
\begin{eqnarray}\label{yfghhgbvhfhggg}
P^{\bot}_{k}u_k\rightharpoonup 0\ \mbox{in}\ X.\end{eqnarray}
Since $X$ can be compactly embedded  into $L^p_r(\mathbb{R}^N)$,
by (\ref{yfghhgbvhfhggg}), we get that for every $2\leq p<2^*,$
\begin{eqnarray}\label{ggfbv88bubufggf}
\lim_{k\rightarrow\infty}\int_{\mathbb{R}^N}|P^\bot_k u_k |^{p}=0.
\end{eqnarray}
By (\ref{hfbv8fyfyyfww}), (\ref{bcbvvvvv55}),
(\ref{ggfbv88bubufggf}) and the condition $\bf (F_1)$, we obtain
\begin{eqnarray}\label{gggvcvbmmiujuj}
\lim_{k\rightarrow\infty}\int_{\mathbb{R}^N}\Big|F(P_k
u_k+\pi_{k}(P_k u_k ))-F(u_k)\Big|=0.\nonumber
\end{eqnarray} Thus
\begin{eqnarray}\label{gvcvbmmiujuj}
\lim_{k\rightarrow\infty}\sup\{\int_{\mathbb{R}^N}\Big|F(P_k
u+\pi_{k}(P_k u ))-F(u)\Big|\ |\
u\in\overline{\mathcal{N}_{\delta,\tau,k}} \}=0.
\end{eqnarray}
By $\bf (ii)$ of Lemma \ref{77eygdjjj876yyyyff5ttt} and
(\ref{gvcvbmmiujuj}), we get that
\begin{eqnarray}\label{jaccxxeedr77}
\lim_{k\rightarrow\infty}||\mathcal{I}_k-I||_{C^0(\overline{\mathcal{N}_{\delta,\tau,k}})}=0.
\end{eqnarray}

For $h\in X,$ \begin{eqnarray}\label{bbv00oqrrw} \langle\nabla
\mathcal{I}_k(u),h\rangle
 &=&\langle u,h\rangle+\langle\pi_{k}(P_k
u), D\pi_{k}(P_k
u)(P_k h)\rangle\nonumber\\
&&-\int_{\mathbb{R}^N}f(P_k u+\pi_{k}(P_k u ))\cdot(P_k h+
D\pi_{k}(P_k u)(P_k h)).\nonumber\end{eqnarray} By $\bf (iii)$ of
Lemma \ref{77eygdjjj876yyyyff5ttt} and the same argument as above,
we can get that
\begin{eqnarray}\label{nbooifgf66t}
\lim_{k\rightarrow\infty}\sup\{ \langle\nabla
\mathcal{I}_k(u)-\nabla I(u),h\rangle\ |\ u\in
\overline{\mathcal{N}_{\delta,\tau,k}},\ || h||\leq 1\}=0.
\end{eqnarray}
The  result of this Lemma follows from (\ref{jaccxxeedr77}) and
(\ref{nbooifgf66t}). \hfill$\Box$

\begin{remark}\label{uty77ty6bcc}  For $r>0,$ let $\sigma\in(0,\sigma_{r/2})$,
where $\sigma_{r/2}$ comes from Lemma \ref{bcng88yuyjhjnew55r},
and let $a\in (c-\sigma,c)$, $b\in(c,c+\sigma)$ be regular values
of $I$, where $c$ comes from (\ref{xiquejiajja}). By Lemma
\ref{bcng88yuyjhjnew55r}, there exists a
 GM pair $(W,W_-)$ of $\mathcal{K}^{b}_{a}$ associated with some
 pseudo-gradient vector field of $I$ such that $W\subset
N_{r/2}(\mathcal{K}^{b}_{a}).$ By (\ref{gdteraa5wdd}), if
$0<r<\min\{\delta,\tau\}$, then
$N_{r}(\mathcal{K})\subset\mathcal{N}_{\delta,\tau,k}$ if $k$ is
sufficiently large. Denote the critical set of $\mathcal{I}_k$ in
$\mathcal{N}_{\delta,\tau,k}$ by $\widehat{\mathcal{K}}_k$. By
$\bf (i)$ and $\bf (iv)$ of Lemma \ref{77eygdjjj876yyyyff5ttt}, we
deduce that $\widehat{\mathcal{K}}_k=P_{k}\mathcal{K}^{b}_{a}$.
Then by (\ref{shengdanweirenjie}),
$\widehat{\mathcal{K}}_{k}\subset \mbox{int}\ W$ if $k$ is large
enough. By \cite[Theorem III.4]{CG} and Lemma \ref{yryuqsaxxx3et},
we infer that for sufficiently large $k$, $(W,W_-)$ is also a GM
pair of $\mathcal{I}_k$ for $\widehat{\mathcal{K}}_k$ associated
with some pseudo-gradient vector filed of $\mathcal{I}_k.$

 For $v\in \mathcal{N}_{\delta,k},$ denote
$I(v+\pi_k(v))$ by $g_k(v).$  And denote the critical  set of
$g_k$ in $W$ by $\mathcal{K}_k$. By $\bf (i)$ and $\bf (iv)$ of
Lemma \ref{77eygdjjj876yyyyff5ttt}, we deduce that
$\mathcal{K}_k=P_{k}\mathcal{K}^{b}_{a}=\widehat{\mathcal{K}}_k$.
Let $(W_k,W^-_k)$ be a GM pair of $g_k$ for $\mathcal{K}_k$. Note
that for $u=w+v\in \mathcal{N}_{\delta,\tau,k}$ with $w\in
X^{\bot}_{k},$ $v\in X_k,$
$\mathcal{I}_k(u)=\frac{1}{2}||w||^2+g_k(v)$. By shifting theorem
(see Lemma 5.1 of \cite{cha}), we have
\begin{eqnarray}\label{000o6fgfga4q3q3gvcbvx}
\check{H}^q(W_k,W^-_k)= \check{H}^q(W,W^-),\
q=0,1,2,\cdots.\nonumber\end{eqnarray} Combining
 Lemma \ref{yryhfgyhgtgftt645t}, we get that, for sufficiently
 large $k,$
\begin{eqnarray}\label{6fgfga4q3q3gvcbvx}
\check{H}^{1}(W_k,W^-_k)= \check{H}^{1}(W,W^-)\neq
0.\end{eqnarray}
\end{remark}

\section{A variational reduction for the  functional $E_\epsilon$}\label{dadczd443eee}
 For $v\in
\cup^{s}_{i=1}B_{X}(u_i,\tau_{u_i})$ and $y\in\mathbb{R}^{N}$,
denote the space
$$\{\zeta(\cdot-y)\ |\ \zeta\in X_{k}\}\oplus \mathcal{T}_v(\cdot-y)$$
  by $T_{v,y,k}$, where $\mathcal{T}_v$ comes from
   (\ref{6erwtsfga55a}). Denote the orthogonal complemental space
of $T_{v,y,k}$ in $Y$  by $T_{v,y,k}^{\bot}$.

Recall that (see (\ref{nnnzbxccc7y}))
\begin{eqnarray}\label{fdt5ermnjhuu}
\mathcal{N}_{\delta,k}=\{u\in X_k\ |\
\mbox{dist}_{X}(u,P_k\mathcal{K})<\delta\}.\nonumber
\end{eqnarray}
For $v\in \mathcal{N}_{\delta,k},$ define
$$L_{v,y,\epsilon,k}: T_{v,y,k}^{\bot}\rightarrow
T_{v,y,k}^{\bot}$$  by
\begin{eqnarray}\label{kfjj7yevcv77y}
  w\in T_{v,y,k}^{\bot}\mapsto w-S_{v,y,k}(-\triangle +1+V(\epsilon
x))^{-1}(f'(v(\cdot-y))w)
\end{eqnarray}
where $S_{v,y,k}: Y\rightarrow T_{v,y,k}^{\bot}$ is orthogonal
projection  and the operator $(-\triangle +1+V(\epsilon x))^{-1}$
is defined by (\ref{uudgttdrhhj77wew}).

\begin{lemma}\label{yhhdvvvcfd4e}
 Given $R > 0,$ there exist  $\delta_0>0$, $\epsilon_0>0$, $l^*>0$ and
$C>0$ which are independent of $k,$ such that   if $k\geq l^*,$
 $0<\delta\leq\delta_0$ and $0\leq\epsilon\leq\epsilon_0,$
then
 for any $v\in \overline{\mathcal{N}_{\delta,k}}$ and $y\in \overline{B_{\mathbb{R}^N}(0,R)}$,
  $L_{v,y,\epsilon,k}$ is
 invertible and
\begin{eqnarray}\label{huchedan8143}||L_{v,y,\epsilon,k}w||\geq C||w||,\ \forall |y|\leq R,\
\forall w\in T_{v,y,k}^{\bot}.\end{eqnarray}
\end{lemma}
\noindent{\bf Proof.} Suppose $\kappa=\max\{\tau_{u_i}\ |\ 1\leq
i\leq s\}$ is small enough such that Lemma \ref{gsfs5s4srfsf}
holds. By (\ref{gdtaazeraa5wdd}), for sufficiently small
$\delta_0>0,$ there exists $l'_\kappa>0$ such that
$\mathcal{N}_{\delta_0,k}\subset\cup^{s}_{i=1}B_X(u_i,\tau_{u_i})$
if $k\geq l'_\kappa.$    Note that $L_{v,0,0,k}$ is exactly the
operator $P_{E^{\bot}_{v,k}}\nabla^2 J(v)|_{E^{\bot}_{v,k}}$ which
has been defined in Lemma \ref{gsfs5s4srfsf} and for every $w\in
T^{\bot}_{v,y,k}$, $$L_{v,y,0,k}w=L_{v,0,0,k}w(\cdot-y).$$ Thus,
by Lemma \ref{gsfs5s4srfsf},
 there exists $C'>0$   such that if
 $k\geq
l^*:=\max\{l_\kappa,l'_\kappa\}$, then for any $v\in
\mathcal{N}_{\delta_0,k}$,
\begin{eqnarray}\label{bbdgvdfdffd88duu}
||L_{v,y,0,k}w||\geq C'||w||,\ \forall |y|\leq R ,\ \forall w\in
T_{v,y,k}^{\bot},\nonumber
\end{eqnarray}
where $l_\kappa$ is the constant comes from Lemma
\ref{gsfs5s4srfsf}. Therefore, to prove (\ref{huchedan8143}), it
suffices to prove that
\begin{eqnarray}\label{bbdgd77dydteef}
&&\lim_{\epsilon\rightarrow
0}\sup\Big\{||L_{v,y,\epsilon,k}w-L_{v,y,0,k}w|| \ |\ w\in
T_{v,y,k}^{\bot},\ ||w||\leq 1,\\
&&\quad\quad\quad\quad\quad
v\in\overline{\mathcal{N}_{\delta_0,k}},\ y\in
\overline{B_{\mathbb{R}^N}(0,R)},\ k\geq l^*\Big\}=0.\nonumber
\end{eqnarray}

If we can prove that for any given sequences
$\{k_n\}\subset\mathbb{N},$ $\{\epsilon_n\}\subset (0,+\infty),$
$\{y_n\}\subset \overline{B_{\mathbb{R}^N}(0,R)}$, $\{v_n\}$ and
$\{w_n\}$ which satisfy that $\epsilon_n\rightarrow 0$ as
$n\rightarrow\infty$, $v_n\in \overline{
\mathcal{N}_{\delta_0,k_{n}}}$,
  $w_n\in T^{\bot}_{v_n,y_n,k_{n}}$ and $||w_n||\leq1$,
  $n=1,2,\cdots$,
\begin{eqnarray}\label{bc8quyuanfuij}
\lim_{n\rightarrow\infty}||L_{v_n,y_n,\epsilon_n,k_{n}}w_n-L_{v_n,y_n,0,k_{n}}w_n||=0,
\end{eqnarray}
then (\ref{bbdgd77dydteef}) holds. We only give the proof of
(\ref{bc8quyuanfuij}) in the case $k_n\rightarrow\infty,$
$n\rightarrow\infty$, since the proofs in  other cases are
similar. Without loss of generality, we assume that $\{k_n\}$ is
exactly the sequence $\{k\}$ and we shall denote $\epsilon_n$,
$y_n,$ $v_n$ and $w_n$ by $\epsilon_k$, $y_k,$ $v_k$ and  $w_k$
respectively, $k=1,2,\cdots.$

Passing to a subsequence, we may assume that as
$k\rightarrow\infty$,   $y_k\rightarrow y_0$, $v_k\rightharpoonup
v_0$ in $X$ and $w_k\rightharpoonup w_0$ in $Y$.

Let $$\eta_k=(-\triangle +1+V(\epsilon_k
x))^{-1}(f'(v_k(\cdot-y_k))w_k).$$ It is easy to verify that
$\{\eta_k\}$ is bounded in $Y$ and
\begin{eqnarray}\label{bbdgcfd66dtdtfef}
\eta_k=(-\triangle+1)^{-1}(f'(v_k(\cdot-y_k))w_k)-
(-\triangle+1)^{-1}V(\epsilon_k)\eta_k.
\end{eqnarray} Passing to a subsequence, we may assume
that $\eta_k\rightharpoonup \eta_0$ in $Y$ as
$k\rightarrow\infty.$

 By  definition of $L_{v,y,\epsilon,k}$ and (\ref{bbdgcfd66dtdtfef}),
  we get that
\begin{eqnarray}\label{ggasd4wewew1111kk}
L_{v_k,y_k,\epsilon,k}w-L_{v_k,y_k,0,k}w=S_{v_k,y_k,k}(-\triangle+1)^{-1}V(\epsilon_k
x)\eta_k.
\end{eqnarray}
The condition $\bf (V_1)$ implies that $V(0)=0$.   It follows that
for any  $h\in Y,$
\begin{eqnarray}\label{ttdgdvc66d5dr}
\lim_{k\rightarrow\infty}\int_{\mathbb{R}^N}V(\epsilon_k x)\eta_k
h=0.\end{eqnarray} Since $\eta_k$ is a weak solution of the
 equation:
\begin{eqnarray}\label{hhdbdf6d5dttd1wwq}
-\triangle\eta_k+\eta_k+V(\epsilon_k
x)\eta_k=f'(v_k(\cdot-y_k))w_k,\end{eqnarray} by
(\ref{ttdgdvc66d5dr}), $y_k\rightarrow y_0,$
$\eta_k\rightharpoonup \eta_0$ and $w_k\rightharpoonup w_0$ in
$Y$, we get that $\eta_0$ is a weak solution of the  equation:
\begin{eqnarray}\label{ggdtdrfdrdrrr}
-\triangle\eta_0+\eta_0=f'(v_0(\cdot-y_0))w_0.\end{eqnarray} From
(\ref{hhdbdf6d5dttd1wwq}) and (\ref{ggdtdrfdrdrrr}), we obtain
\begin{eqnarray}\label{gdvct55drsss}
&&-\triangle(\eta_k-\eta_0)+(\eta_k-\eta_0)+V(\epsilon_k
x)(\eta_k-\eta_0)\nonumber\\
&=&(f'(v_k(\cdot-y_k))w_k-f'(v_0(\cdot-y_0))w_0)-V(\epsilon_k
x)\eta_0.\nonumber\end{eqnarray} Multiplying the above equation by
$\eta_k-\eta_0$ and integrating, we get that there exists a
constant $C>0$ such that
\begin{eqnarray}\label{fad4adad77auueear}
&&C||\eta_k-\eta_0||^2\nonumber\\
&\leq&||\eta_k-\eta_0||^2+\int_{\mathbb{R}^N}V(\epsilon_k
x)(\eta_k-\eta_0)^2\ (\mbox{by the condition }\ {\bf (V_0)})\nonumber\\
&=&\int_{\mathbb{R}^N}\Big(f'(v_k(\cdot-y_k))w_k-f'(v_0(\cdot-y_0))w_0-V(\epsilon_k
x)\eta_0\Big)\cdot(\eta_k-\eta_0)\nonumber\\
&\leq&\int_{\mathbb{R}^N}\Big|f'(v_k(\cdot-y_k))w_k-f'(v_0(\cdot-y_0))w_0\Big|\cdot
|\eta_k-\eta_0|\nonumber\\
&&+(\int_{\mathbb{R}^N}V^2(\epsilon_k
x)\eta^{2}_{0})^{\frac{1}{2}}\cdot
||\eta_k-\eta_0||_{L^{2}(\mathbb{R}^N)}.
\end{eqnarray}
Since $v_k\rightharpoonup v_0$ in $X$ and $y_k\rightarrow y_0$ as
$k\rightarrow\infty$, by the fact that $X$ can be compactly
embedding
 into $L^{p}_{r}(\mathbb{R}^N)$ ($\forall p\in [2,2^*)$),
we get that
\begin{eqnarray}\label{shengshengbuxiew4}
\lim_{k\rightarrow\infty}||v_k(\cdot-y_k)-v_0(\cdot-y_0)||_{L^p(\mathbb{R}^N)}=0,\
\forall p\in [2,2^*).
\end{eqnarray}
By (\ref{shengshengbuxiew4}) and the condition $\bf (F_1)$, we get
that
\begin{eqnarray}\label{und88qypdmchhhg555}
\lim_{k\rightarrow\infty}\int_{\mathbb{R}^N}\Big|f'(v_k(\cdot-y_k))w_k-f'(v_0(\cdot-y_0))w_0\Big|\cdot
|\eta_k-\eta_0|=0.
\end{eqnarray}
By (\ref{fad4adad77auueear}), (\ref{und88qypdmchhhg555}) and
\begin{eqnarray}\label{ggdbc66s5srrrs99i}
\lim_{k\rightarrow\infty}\int_{\mathbb{R}^N}V^2(\epsilon_k
x)\eta^{2}_{0}=0,\end{eqnarray} we get that
\begin{equation}\label{bd7dgsfsluoo01wwq}
\lim_{k\rightarrow\infty}||\eta_k-\eta_0||=0.
\end{equation}
(\ref{ggdbc66s5srrrs99i}) and (\ref{bd7dgsfsluoo01wwq}) yield
\begin{eqnarray}\label{kkkkkkkata554}
\lim_{k\rightarrow\infty}\int_{\mathbb{R}^N}V^2(\epsilon_k
x)\eta^{2}_{k}=0.\end{eqnarray} It follows that
\begin{eqnarray}\label{hdudtekmfhfg}
\lim_{k\rightarrow\infty}||(-\triangle+1)^{-1}V(\epsilon_k
x)\eta_{k}||=0.
\end{eqnarray}
Combining (\ref{hdudtekmfhfg}) and (\ref{ggasd4wewew1111kk}) leads
to (\ref{bc8quyuanfuij}).

\medskip

Finally, by definition,   $L_{v,y,\epsilon,k}$ is a Fredholm
operator with index zero and by (\ref{huchedan8143}),  it is an
injection. Therefore, it is invertible. \hfill$\Box$

\begin{Theorem}\label{77tutyyy} Given $R>0.$ There exist
$\delta^*>0$ and $\epsilon^*>0$ such that if
$0<\delta\leq\delta^*$ and $0\leq\epsilon\leq\epsilon^*$, then
there exist $k(\delta)$ and a $C^{1}-$mapping
$$w_{\delta,k}(\cdot,\cdot,\epsilon):\overline{\mathcal{N}_{\delta, k}}
\times\overline{B_{\mathbb{R}^N}(0,R)} \rightarrow  Y, \
(u,y)\mapsto w_{\delta,k}(u,y,\epsilon)$$ for $k\geq k(\delta)$,
 satisfying
\begin{itemize}

\item[{\bf (i)}] $w_{\delta,k}(u,y,\epsilon)\in T_{u,y,k}^{\bot},$
$\forall (u,y)\in \overline{\mathcal{N}_{\delta, k}}\times
\overline{B_{\mathbb{R}^{N}}(0,R)};$

\item[{\bf (ii)}] $\langle \nabla
E_{\epsilon}(u(\cdot-y)+w_{\delta,k}(u,y,\epsilon)),
\phi\rangle=0,$ $\forall \phi\in T_{u,y,k}^{\bot};$

\item[{\bf (iii)}] $w_{\delta,k}(u,y,0)=(\pi_{k}(u))(\cdot-y),$
$\forall (u,y)\in \overline{\mathcal{N}_{\delta, k}}\times
\overline{B_{\mathbb{R}^{N}}(0,R)}$;

\item[{\bf (iv)}] for any $r>0,$ there exists $\delta_r>0$ such
that if $0<\delta\leq\delta_r$,
$u\in\overline{\mathcal{N}_{\delta,k}}$,
$y\in\overline{B_{\mathbb{R}^N}(0,R)}$ and $k\geq k(\delta)$, then
$||w_{\delta,k}(u,y,\epsilon)||\leq r;$

\item[{\bf (v)}]for any    $n>0$,
\begin{eqnarray}\label{gdbcgtd55drsf}
&&\sup\{
 ||(1+|x|)^{n}w_{\delta,k}(u,y,\epsilon)||_{L^{\infty}(\mathbb{R}^N)}\ |\
 (u,y)\in\overline{\mathcal{N}_{\delta,k}}\times\overline{B_{\mathbb{R}^N}(0,R)},\
0\leq\epsilon\leq\epsilon^*\}\nonumber\\
&&<\infty.
\end{eqnarray}
\end{itemize}
\end{Theorem}
\noindent{\bf Proof.} By Lemma  \ref{yhhdvvvcfd4e}, we know that
for any $R>0$,
   $L_{u,y,\epsilon,k}$ is invertible
   if  $0<\delta\leq\delta_0$, $0\leq\epsilon\leq\epsilon_0$ and  $k\geq l^*$. Moreover,
    the upper bound of $
||L^{-1}_{u,y,\epsilon,k}||$ is independent of $u,$ $y$,
$\epsilon$ and $k.$   For $u\in\overline{\mathcal{N}_{\delta,k}}$
and $r>0$, let
$$\Phi_{u,y,\epsilon,k}:\overline{B_{T_{u,y,k}^{\bot}}(0,r)}\rightarrow T_{u,y,k}^{\bot},$$
$$w\mapsto w-L^{-1}_{u,y,\epsilon,k}S_{u,y,k}\nabla E_{\epsilon}(u(\cdot-y)+w).$$

Now, we show that if $r$, $\delta$ and  $\epsilon$ are
 small enough and $k$ is large enough, then for any
$u\in\overline{\mathcal{N}_{\delta,k}}$,
 $\Phi_{u,y,\epsilon,k}$ is a contractive mapping in $\overline{B_{T_{u,y,k}^{\bot}}(0,r)}.$

Using
\begin{eqnarray}\label{uhdvcvv66t1q}
&&\nabla E_{\epsilon}(u(\cdot-y)+w)\nonumber\\
&=&u(\cdot-y)+w-(-\triangle +1+V(\epsilon
x))^{-1}f(u(\cdot-y)+w)\nonumber
\end{eqnarray}
and the mean value theorem, we get that for any $w_1,
w_2\in\overline{B_{T_{u,y,k}^{\bot}}(0,r)}$,
$\Phi_{u,y,\epsilon,k}(w_1)-\Phi_{u,y,\epsilon,k}(w_2)$ equals
\begin{eqnarray}\label{bbdvcfdfhhadda77253rr}
&&(w_1-w_2)-L^{-1}_{u,y,\epsilon,k}S_{u,y,k}
\Big\{(w_1-w_2)\nonumber\\&&-(-\triangle +1+V(\epsilon
x))^{-1}(f'(u(\cdot-y)+\tilde{w})\cdot(w_1-w_2))\Big\}\nonumber\\
&=&(w_1-w_2)-L^{-1}_{u,y,\epsilon,k}S_{u,y,k}
\Big\{(w_1-w_2)\nonumber\\&&-(-\triangle +1+V(\epsilon
x))^{-1}f'(u(\cdot-y))(w_1-w_2)\nonumber\\
&&-(-\triangle +1+V(\epsilon
x))^{-1}(f'(u(\cdot-y)+\tilde{w})-f'(u(\cdot-y)))(w_1-w_2)\Big\}
\end{eqnarray}
where $\tilde{w}=\theta w_1+(1-\theta)w_2$ for some $0<\theta<1.$
By the condition $\bf (F_1)$, we can prove that
\begin{eqnarray}\label{ccsdseee6t6t}
&&\lim_{r\rightarrow 0}\sup\{ ||(-\triangle +1+V(\epsilon
x))^{-1}(f'(u(\cdot-y)+\tilde{w})-f'(u(\cdot-y)))\varphi||\\
&&\quad\quad\quad\quad\ |\ u\in\overline{\mathcal{N}_{\delta,k}},\
|y|\leq R,\ \varphi\in Y,\ ||\varphi||\leq 1,\
0\leq\epsilon\leq\epsilon_0\}=0.\nonumber
\end{eqnarray}
 By
$||L^{-1}_{u,y,\epsilon,k}||_{\mathcal{L}(Y)}\leq 1/C$ (see Lemma
\ref{yhhdvvvcfd4e} ), $||S_{u,y,k}|| _{\mathcal{L}(Y)}\leq 1$ and
(\ref{ccsdseee6t6t}), we deduce that if $r$  is small enough, then
\begin{eqnarray}\label{ccsdseaas11ee6t6tjjafaf}
&& ||L^{-1}_{u,y,\epsilon,k}S_{u,y,k}(-\triangle +1+V(\epsilon
x))^{-1}(f'(u(\cdot-y)+\tilde{w})-f'(u(\cdot-y)))(w_1-w_2)||\\
&&\leq\frac{1}{C} ||(-\triangle +1+V(\epsilon
x))^{-1}(f'(u(\cdot-y)+\tilde{w})-f'(u(\cdot-y)))(w_1-w_2)||\nonumber\\
&&\leq\frac{1}{2}||w_1-w_2||.\nonumber\end{eqnarray} By the
definition of $L_{u,y,\epsilon,k}$,
\begin{eqnarray}\label{bbdvcfjjjjdf77253rr}
&&L^{-1}_{u,y,\epsilon,k}S_{u,y,k} \Big\{(w_1-w_2)-(-\triangle
+1+V(\epsilon
x))^{-1}(f'(u(\cdot-y))(w_1-w_2))\Big\}\\
&&=(w_1-w_2).\nonumber
\end{eqnarray}
Combining (\ref{ccsdseaas11ee6t6tjjafaf}),
(\ref{bbdvcfjjjjdf77253rr}) and (\ref{bbdvcfdfhhadda77253rr}), we
deduce that there exists $r_0>0$ such that if  $0<r\leq r_0$,
$0<\delta\leq\delta_0$, $0\leq\epsilon\leq\epsilon_0$ and $k\geq
l^*$, then for any
$(u,y)\in\overline{\mathcal{N}_{\delta,k}}\times
\overline{B_{\mathbb{R}^N}(0,R)}$ and $w_1, w_2\in
\overline{B_{T_{u,y,k}^{\bot}}(0,r)},$
\begin{equation}\label{taiushdgvv}
||\Phi_{u,y,\epsilon,k}(w_1)-\Phi_{u,y,\epsilon,k}(w_2)||\leq\frac{1}{2}
||w_1-w_2||.\end{equation}

\noindent{\bf Claim:} For any $0<r\leq r_0,$ there exist
$\epsilon_r$, $\delta_r$ and $k(\delta,r)$  such that if
$0<\delta\leq\delta_r$, $0\leq\epsilon\leq\epsilon_r$ and $k\geq
k(\delta,r)$, then
\begin{eqnarray}\label{qquauagt55e4881i}
||\Phi_{u,y,\epsilon,k}(0)||\leq r/2,\  \forall
(u,y)\in\overline{\mathcal{N}_{\delta,k}}\times\overline{B_{\mathbb{R}^N}(0,R)}.
\end{eqnarray}

 Let $h_{u,y,\epsilon}=(-\triangle
+1+V(\epsilon x))^{-1}f(u(\cdot-y)).$ It is easy to verify
\begin{eqnarray}\label{nnfb99273654} h_{u,y,\epsilon}=
(-\triangle +1)^{-1}f(u(\cdot-y))-(-\triangle +1)^{-1}V(\epsilon
x)h_{u,y,\epsilon}.
\end{eqnarray}
The same argument as (\ref{kkkkkkkata554}) yields
$$\lim_{\epsilon\rightarrow 0}\sup\{\int_{\mathbb{R}^N}V^2(\epsilon x)h^2_{u,y,\epsilon}\ |\
\ u\in \overline{\mathcal{N}_{\delta_{_{0}},k}},\ y\in\overline{
B_{\mathbb{R}^N}(0,R)},\ k\geq l^*\}=0.$$  Thus, by
(\ref{nnfb99273654}), as $\epsilon\rightarrow 0,$
\begin{eqnarray}\label{ffdvcgttr5553red}
&&\sup\{||(-\triangle +1+V(\epsilon
x))^{-1}f(u(\cdot-y))\nonumber\\
&&\quad\quad\quad -(-\triangle +1)^{-1}f(u(\cdot-y))||\ |\ u\in
\overline{\mathcal{N}_{\delta_{_{0}},k}},\ y\in\overline{
B_{\mathbb{R}^N}(0,R)},\ k\geq l^*\}\nonumber\\
&&\rightarrow 0.\nonumber
\end{eqnarray}
It follows that as $\epsilon\rightarrow 0,$
\begin{eqnarray}\label{ffdaqqqqqed444453red}
&&\sup\{||\nabla E_\epsilon(u(\cdot-y))-\nabla J(u(\cdot-y))||\ |\
u\in \overline{\mathcal{N}_{\delta_{_{0}},k}},\ y\in\overline{
B_{\mathbb{R}^N}(0,R)},\ k\geq l^*\}\\
&&\rightarrow 0.\nonumber
\end{eqnarray}
Therefore, for $0<r\leq r_0$, there exists $\epsilon_r>0$ such
that for any $u\in \overline{\mathcal{N}_{\delta_{_{0}},k}},$ $
y\in\overline{ B_{\mathbb{R}^N}(0,R)}$ and $ k\geq l^*$,
\begin{eqnarray}\label{tgf77fyr6tzcmv}
||\nabla E_\epsilon(u(\cdot-y))-\nabla
J(u(\cdot-y))||<\frac{C}{4}r\ \ \mbox{if}\
0\leq\epsilon\leq\epsilon_r,
\end{eqnarray}
where the constant $C$ comes from Lemma \ref{yhhdvvvcfd4e}. Since
$\nabla J(v(\cdot-y))=\nabla J(v)=0,$ $\forall v\in\mathcal{K},$
we get that for any $0<r\leq r_0$, there exists $\delta_r$ such
that for any $0<\delta\leq\delta_r$ and any $u\in
N_{2\delta}(\mathcal{K})$,
\begin{eqnarray}\label{112qtgf7jjxiangjian} ||\nabla
J(u(\cdot-y))||<\frac{C}{4}r.
\end{eqnarray}By (\ref{112qtgf7jjxiangjian}) and
 the fact that (see (\ref{shengdanweirenjie}))
$$\lim_{k\rightarrow\infty}\overline{\mathcal{N}_{\delta,k}}\subset
N_{2\delta}(\mathcal{K}),$$ we deduce that there exists
$k(\delta,r)$ such that if $k\geq k(\delta,r)$, then for any
$0<\delta\leq\delta_r$ and any $u\in
\overline{\mathcal{N}_{\delta,k}}$,
\begin{eqnarray}\label{112qtgf77fyr6tzcmv} ||\nabla
J(u(\cdot-y))||<\frac{C}{4}r.
\end{eqnarray}
Thus, the claim follows from (\ref{tgf77fyr6tzcmv}),
(\ref{112qtgf77fyr6tzcmv}) and the fact that
$$||\Phi_{u,y,\epsilon,k}(0)||\leq \frac{1}{C}||\nabla
E_\epsilon(u(\cdot-y))||.$$ Combining (\ref{taiushdgvv}) and
(\ref{qquauagt55e4881i}) leads to
$$||\Phi_{u,y,\epsilon,k}(w)||\leq r$$ for every $w\in
\overline{B_{T_{u,y,k}^{\bot}}(0,r)}$. Therefore,
$\Phi_{u,y,\epsilon,k}$ is a contractive mapping in
$\overline{B_{T_{u,y,k}^{\bot}}(0,r)}.$
 By Banach fixed point theorem, there exists unique fixed point
$w_{\delta,k}(u,y,\epsilon)$ of $\Phi_{u,y,\epsilon,k}$. Denote
$\delta_{r_{_{0}}}$ by $\delta^*$, $\epsilon_{r_{_{0}}}$ by
$\epsilon^*$ and $k(\delta, r_0)$ by $k(\delta)$. It is easy to
verify that the conclusions $\bf (i)-(iv)$ hold for
$w_{\delta,k}(u,y,\epsilon)$.

Now, we   prove that
$w_{\delta,k}:\overline{\mathcal{N}_{\delta,k}}\times
\overline{B_{\mathbb{R}^{N}}(0,R) }\rightarrow  Y$ is $C^1$. For
any $(u_0,y_0)\in \overline{\mathcal{N}_{\delta,k}}\times
\overline{B_{\mathbb{R}^{N}}(0,R)}$ and $(u,y)$ close to
$(u_0,y_0)$, both
$S_{u_0,y_0,k}|_{T^\perp_{u,y,k}}:T^\perp_{u,y,k}\to
T^\perp_{u_0,y_0,k}$ and
$S_{u,y,k}|_{T^\perp_{u_0,y_0,k}}:T^\perp_{u_0,y_0,k}\to
T^\perp_{u,y,k}$ are isomorphisms, and finding a solution $w\in
T^{\perp}_{u,y,k}$ to the equation $S_{u,y,k}\nabla
E_{\epsilon}(u(\cdot-y)+w)=0$ is equivalent to finding a solution
$w\in T^{\bot}_{u_0,y_0,k}$ to the equation
$S_{u_{0},y_0,k}S_{u,y,k}\nabla
E_{\epsilon}(u(\cdot-y)+S_{u,y,k}w)=0$. Note that
$S_{u_0,y_0,k}S_{u,y,k}\nabla E_{\epsilon}(u(\cdot-y)+S_{u,y,k}w)$
is $C^1$ near $(u_0, y_0,w_0)\in
\overline{\mathcal{N}_{\delta,k}}\times
\overline{B_{\mathbb{R}^{N}}(0,R)}\times T^{\perp}_{u_0,y_0,k}$
and the Fr\'echet partial derivative of
$S_{u_0,y_0,k}S_{u,y,k}\nabla E_{\epsilon}(u(\cdot-y)+S_{u,y,k}w)$
at $(u_0,y_0,w_0)$ with respect to $w$ is $L_{u_0,y_0,\epsilon,k}$
which is invertible. Therefore, the implicit functional theorem
implies that
$$w_{\delta,k}(\cdot,\cdot,\epsilon):\overline{\mathcal{N}_{\delta,k}}\times
\overline{B_{\mathbb{R}^{N}}(0,R)} \rightarrow  Y$$ is $C^1$.

Finally, we give the  proof of $\bf (v).$ Let
\begin{eqnarray}\label{agbsys67stfsdfuh}
\varphi_{u,y,\epsilon,k}=u(\cdot-y)+w_{\delta,k}(u,y,\epsilon)-
P_{T_{u,y,k}}(\nabla
E_\epsilon(u(\cdot-y)+w_{\delta,k}(u,y,\epsilon))),
\end{eqnarray}
where $P_{T_{u,y,k}}:Y\rightarrow T_{u,y,k}$ is  orthogonal
projection. By the conclusion $\bf (ii)$ of this Theorem, we get
that
\begin{eqnarray}\label{udg55dreweew55uh}
P_{T_{u,y,k}}(\nabla
E_\epsilon(u(\cdot-y)+w_{\delta,k}(u,y,\epsilon)))=\nabla
E_\epsilon(u(\cdot-y)+w_{\delta,k}(u,y,\epsilon)).\end{eqnarray}
Thus, by (\ref{agbsys67stfsdfuh}) and (\ref{udg55dreweew55uh}),
$\varphi_{u,y,\epsilon,k}$ satisfies
\begin{eqnarray}\label{vdvgftuettet}
&&-\triangle\varphi_{u,y,\epsilon,k}+\varphi_{u,y,\epsilon,k}+V(\epsilon
x)\varphi_{u,y,\epsilon,k}=f(u(\cdot-y)+w_{\delta,k}(u,y,\epsilon)).
\end{eqnarray}
By the  definition of $ T_{u,y,k}$, we have
\begin{eqnarray}\label{bbbc77yhdtdfffd}
&&P_{T_{u,y,k}}(\nabla
E_\epsilon(u(\cdot-y)+w_{\delta,k}(u,y,\epsilon)))\nonumber\\
&=&\sum^{N}_{j=1}\Big\langle\nabla
E_\epsilon(u(\cdot-y)+w_{\delta,k}(u,y,\epsilon)),\sum^{s}_{i=1}\xi_i(u)\frac{u_i(\cdot-y)}{\partial
x_j}\Big\rangle\frac{\sum^{s}_{i=1}\xi_i(u)\frac{u_i(\cdot-y)}{\partial
x_j}}{||\sum^{s}_{i=1}\xi_i(u)\frac{u_i(\cdot-y)}{\partial
x_j}||^2}\nonumber\\
&&+\sum^{k}_{i=1}\langle\nabla
E_\epsilon(u(\cdot-y)+w_{\delta,k}(u,y,\epsilon)),\tilde{e}_{i,k}(\cdot-y)\rangle
\tilde{e}_{i,k}(\cdot-y)\nonumber\\
&&+\sum^{q}_{i=1}\langle\nabla
E_\epsilon(u(\cdot-y)+w_{\delta,k}(u,y,\epsilon)),e_i(\cdot-y)\rangle
e_i(\cdot-y).
\end{eqnarray}
Since $\tilde{e}_{i,k}$, $e_i$, $u$  and $\frac{\partial
u_i}{\partial x_j}$ satisfy exponential decay  at infinity, by
(\ref{bbbc77yhdtdfffd}),  for any given $k\geq k(\delta)$ and
$n\geq0$, there exists $C'_{n,k}>0$ such that
\begin{eqnarray}\label{gdbyqytwr00of}
&&\sup\{||(1+|x|)^{n}(P_{T_{u,y,k}}(\nabla
E_\epsilon(u(\cdot-y)+w_{\delta,k}(u,y,\epsilon))))||_{L^{\infty}(\mathbb{R}^N)}\nonumber\\
&&\quad\quad |\
u\in\overline{\mathcal{N}_{\delta,k}},y\in\overline{B_{\mathbb{R}^N}(0,R)},0\leq\epsilon\leq\epsilon^*\}
 \leq
 C'_{k,n}
\end{eqnarray}
and \begin{eqnarray}\label{gdbcgtdrrq55drsf}
\sup_{u\in\overline{\mathcal{N}_{\delta,k}},y\in\overline{B_{\mathbb{R}^N}(0,R)}
 }||(1+|x|)^{n}u(\cdot-y)||_{L^{\infty}(\mathbb{R}^N)}\leq
 C'_{k,n}.
\end{eqnarray}
 Note that $\varphi_{u,y,\epsilon,k}$ satisfies the elliptic equation
(\ref{vdvgftuettet}). Therefore, by the bootstrap argument  and
the fact that
$$\{w_{\delta,k}(u,y,\epsilon)) \ |\
u\in\overline{\mathcal{N}_{\delta,k}},\
y\in\overline{B_{\mathbb{R}^N}(0,R)},\
0\leq\epsilon\leq\epsilon^*\}$$ is compact in $Y$ (because for
fixed $k$, $\overline{\mathcal{N}_{\delta,k}}$ is compact), we get
that
\begin{eqnarray}\label{meiyohb666uquanjjus77ys}
\sup\{||\varphi_{u,y,\epsilon,k}||_{L^{\infty}(\mathbb{R}^N)}\ |\
u\in \overline{\mathcal{N}_{\delta,k}},\
y\in\overline{B_{\mathbb{R}^N}(0,R)},\
0\leq\epsilon\leq\epsilon^*\}<\infty
\end{eqnarray}
and
\begin{eqnarray}\label{bbd3gdfttdrdrrr}
\lim_{\rho\rightarrow\infty}\sup\{||\varphi_{u,y,\epsilon,k}||_{L^{\infty}(\mathbb{R}^N\setminus
\overline{B_{\mathbb{R}^N}(0,\rho)})}\ |\ u\in
\overline{\mathcal{N}_{\delta,k}},\
y\in\overline{B_{\mathbb{R}^N}(0,R)},\
0\leq\epsilon\leq\epsilon^*\}=0.
\end{eqnarray}
By (\ref{meiyohb666uquanjjus77ys}), (\ref{bbd3gdfttdrdrrr}) and
(\ref{agbsys67stfsdfuh}), we get that
\begin{eqnarray}\label{bbdgdfhhhhttdrdrrr}
&&\sup\{||w_{\delta,k}(u,y,\epsilon)||_{L^{\infty}(\mathbb{R}^N)}\
|\ u\in \overline{\mathcal{N}_{\delta,k}},\
y\in\overline{B_{\mathbb{R}^N}(0,R)},\
0\leq\epsilon\leq\epsilon^*\}<\infty.
\end{eqnarray}
and
\begin{eqnarray}\label{bbdgdfttdrdrrr}
&&\lim_{\rho\rightarrow\infty}\sup\{||w_{\delta,k}(u,y,\epsilon)||_{L^{\infty}(\mathbb{R}^N\setminus
\overline{B_{\mathbb{R}^N}(0,\rho)})}\ |\ u\in
\overline{\mathcal{N}_{\delta,k}},\
y\in\overline{B_{\mathbb{R}^N}(0,R)},\
0\leq\epsilon\leq\epsilon^*\}\\ &&=0.\nonumber
\end{eqnarray}
Let $d(t)=f(t)/t,$ $t\in\mathbb{R}.$ Then by
(\ref{bbdgdfhhhhttdrdrrr}), (\ref{gdbcgtdrrq55drsf}) and the
condition $\bf (F_1)$, we have
\begin{eqnarray}\label{fc33fhhhhttdrdrrr}
&&\sup\{||d(u(\cdot-y)+w_{\delta,k}(u,y,\epsilon))||_{L^{\infty}(\mathbb{R}^N)}\
|\ u\in \overline{\mathcal{N}_{\delta,k}},\
y\in\overline{B_{\mathbb{R}^N}(0,R)},\
0\leq\epsilon\leq\epsilon^*\}\\
&&<\infty.\nonumber
\end{eqnarray}
By the condition $\bf (V_0)$, the condition $\bf (F_1)$ and
(\ref{bbdgdfttdrdrrr}), we deduce that there exists $\rho_0$ such
that
\begin{eqnarray}\label{bbdgdt66tssdddd}
&&\inf\{1+V(\epsilon x)-d(u(x-y)+w_{\delta,k}(u,y,\epsilon))\
|\ |x|>\rho_0,\ u\in \overline{\mathcal{N}_{\delta,k}},\nonumber\\
&&\quad\quad\ y\in\overline{B_{\mathbb{R}^N}(0,R)},\
0\leq\epsilon\leq\epsilon^*\}>0.
\end{eqnarray}
Let $\eta$ be a cut-off function which satisfies that $\eta\equiv
1$ in $B_{\mathbb{R}^N}(0,\rho_0)$ and $\eta\equiv 0$ in
$\mathbb{R}^N\setminus \overline{B_{\mathbb{R}^N}(0,\rho_0+1)}$.
We can rewrite equation (\ref{vdvgftuettet}) as
\begin{eqnarray}\label{bbdgcvffdrdrs}
&&-\triangle\varphi_{u,y,\epsilon,k}+(1+V(\epsilon
x)-(1-\eta(x))d(u(x-y)+w_{\delta,k}(u,y,\epsilon)))\varphi_{u,y,\epsilon,k}\\
&=&f_{u,y,\epsilon,k}\nonumber
\end{eqnarray}
with
\begin{eqnarray}\label{hhdnggggasasa}
f_{u,y,\epsilon,k}
&=&d(u(\cdot-y)+w_{\delta,k}(u,y,\epsilon))\cdot u(\cdot-y)\nonumber\\
&&
+\eta(x)\cdot d(u(\cdot-y)+w_{\delta,k}(u,y,\epsilon))\cdot w_{\delta,k}(u,y,\epsilon)\nonumber\\
&&-(1-\eta(x))\cdot
d(u(\cdot-y)+w_{\delta,k}(u,y,\epsilon))\nonumber\\
&&\quad\quad\quad\times(u(\cdot-y)-P_{T_{u,y,k}}(\nabla
E_\epsilon(u(\cdot-y)+w_{\delta,k}(u,y,\epsilon))).
\end{eqnarray}
By (\ref{gdbcgtdrrq55drsf}), (\ref{gdbyqytwr00of}),
(\ref{fc33fhhhhttdrdrrr})  and the fact that
$$\eta(x)d(u(\cdot-y)+w_{\delta,k}(u,y,\epsilon))\cdot w_{\delta,k}(u,y,\epsilon)$$
has compact support, we deduce that there exists $C'''_{n,k}>0$
such that \begin{eqnarray}\label{qqqqg7dbcgtdrrq55drsf}
\sup_{u\in\overline{\mathcal{N}_{\delta,k}},y\in\overline{B_{\mathbb{R}^N}(0,R)}
 }||(1+|x|)^{n}f_{u,y,\epsilon,k}||_{L^{\infty}(\mathbb{R}^N)}\leq
 C'''_{k,n}.
\end{eqnarray}
By (\ref{qqqqg7dbcgtdrrq55drsf}), (\ref{bbdgdt66tssdddd}),
(\ref{bbdgcvffdrdrs})
 and
\cite[Proposition 4.2]{zhang qi s}, we get that there exists
$C''_{n,k}>0$ such that
\begin{eqnarray}\label{qqqqgdbcgtdrrq55drsf}
\sup_{u\in\overline{\mathcal{N}_{\delta,k}},y\in\overline{B_{\mathbb{R}^N}(0,R)}
 }||(1+|x|)^{n}\varphi_{u,y,\epsilon,k}||_{L^{\infty}(\mathbb{R}^N)}\leq
 C''_{k,n}.
\end{eqnarray}
Then the  conclusion  $\bf (v)$ follows from
(\ref{agbsys67stfsdfuh}), (\ref{qqqqgdbcgtdrrq55drsf}),
 (\ref{gdbyqytwr00of}) and (\ref{gdbcgtdrrq55drsf}).
\hfill$\Box$

\medskip

By  the conclusion $\bf (iii)$ of Theorem  \ref{77tutyyy}, we get
that
\begin{eqnarray}\label{bbcfdrref77ayar}
J(u(\cdot-y)+w_{\delta,k}(u,y,0))\equiv I(u+\pi_{k}(u)),\ \forall
(u,y)\in\overline{\mathcal{N}_{\delta,
k}}\times\overline{B_{\mathbb{R}^N}(0,R)}.\end{eqnarray}

In what follows, for a $C^1$ mapping $f$ defined in
$\mathcal{N}_{\delta,k}\times B_{\mathbb{R}^N}(0,R)$, we use the
the notations     $Df$, $D_u f$ and $D_y f$ to denote the
derivatives of $f$  with respect to $(u,y)$ variable, $u$ variable
and $y$ variable respectively and use $Df(u,y)[\bar{u},\bar{y}]$
to denote the derivative of $f$ at the point $(u,y)$ along the
vector $(\bar{u},\bar{y})\in X_k\times \mathbb{R}^N.$ Furthermore,
we use  $D_u f(u,y)[\bar{u}]$ and $D_yf(u,y)[\bar{y}]$ to denote
the Fr\'echet partial derivatives with respect to the $u$ and $y$
variables along the vectors $\bar{u}$ and $\bar{y}$ respectively.

  The condition $\bf (V_1)$ for the potential $V$ yields
\begin{eqnarray}\label{nvb77dydg8uidjnquzhi}
\lim_{\epsilon\rightarrow 0}\frac{V(\epsilon
x)}{\epsilon^{n^*}}=Q_{n^*}(x).\end{eqnarray} The proof of the
following proposition will be  given in appendix.
\begin{proposition}\label{gdgvcfd6d5dttdt}
Let $\delta>0$ be sufficiently small and $k\geq k(\delta)$. If
$\iota<n^*,$ then
\begin{eqnarray}\label{gwts55sr66qeadd}
&&\lim_{\epsilon\rightarrow
0}\sup\{\frac{1}{\epsilon^\iota}\Lambda_k(u,y,\epsilon)\ |\
(u,y)\in\overline{\mathcal{N}_{\delta,
k}}\times\overline{B_{\mathbb{R}^{N}}(0,R)}\}=0\nonumber
\end{eqnarray}
where \begin{eqnarray}\label{hhhfgbfddsjianshanri}
\Lambda_k(u,y,\epsilon)
&=&||w_{\delta,k}(u,y,\epsilon)-\pi_k(u)(\cdot-y)||\nonumber\\
&& +\sup_{\bar{y}\in \mathbb{R}^N,|\bar{y}|\leq
1}||Dw_{\delta,k}(u,y,\epsilon)[0,\bar{y}]-D(\pi_k(u)(\cdot-y))[0,\bar{y}]||\nonumber\\
&&+\sup_{v\in X_k,||v||\leq
1}||Dw_{\delta,k}(u,y,\epsilon)[v,0]-D(\pi_k(u)(\cdot-y))[v,0]||.\nonumber\end{eqnarray}
 Moreover, there
exists a constant $M>0$ which is independent of $(u,y)$ and
$\epsilon$ such that for every
$(u,y)\in\overline{\mathcal{N}_{\delta,
k}}\times\overline{B_{\mathbb{R}^{N}}(0,R)}$ and
$0\leq\epsilon\leq\epsilon^*,$
\begin{eqnarray}\label{gwts55sr66qeaddqqqw}
\Lambda_k(u,y,\epsilon)\leq
M\epsilon^{n^*}.\nonumber\end{eqnarray}
\end{proposition}

 For $0<\delta\leq\delta^*$ and $0\leq\epsilon\leq\epsilon^*,$
 denote the functional \begin{equation}\label{hygggvxv555qmmbn}
E_{\epsilon}(u(\cdot-y)+w_{\delta,k}(u,y,\epsilon)),\ (u,y)\in
\overline{\mathcal{N}_{\delta, k}}\times
\overline{B_{\mathbb{R}^{N}}(0,R)}\end{equation} by
$\Psi_k(u,y,\epsilon)$.

\begin{Theorem}\label{ttrgfg77y} Suppose that  $0<\delta\leq\delta^*$ and $k\geq
k(\delta)$. Then there exists $\epsilon_k>0$ such that if
$0\leq\epsilon\leq\epsilon_k$ and $(u_{\epsilon},y_{\epsilon})\in
\mathcal{N}_{\delta, k}\times B_{\mathbb{R}^{N}}(0,R)$ is a
critical point of the functional $\Psi_k(u,y,\epsilon)$, that is,
\begin{eqnarray}\label{hronrenhfbvgyf6rttrf}
D\Psi_k(u_\epsilon,y_\epsilon,\epsilon)[v,\bar{y}]=0,\ \forall
(v,\bar{y})\in X_k\times \mathbb{R}^N,
\end{eqnarray}
then
$u_{\epsilon}(\cdot-y_{\epsilon})+w_{\delta,k}(u_{\epsilon},y_{\epsilon},\epsilon)$
is a critical point of $E_{\epsilon}$.
\end{Theorem}
\noindent{\bf Proof.} By  the conclusion $\bf (ii)$ of Theorem
\ref{77tutyyy} and hypothesis (\ref{hronrenhfbvgyf6rttrf}),
 we deduce that to prove
 $u_{\epsilon}(\cdot-y_{\epsilon})+w_{\delta,k}(u_{\epsilon},y_{\epsilon},\epsilon)$
 is a critical point of $E_{\epsilon},$ it suffices to prove
 that for sufficiently small $\epsilon>0,$
 \begin{eqnarray}\label{vvcbyryryfg6}
&&\{v(\cdot-y_\epsilon)-(\bar{y}\cdot\nabla_x
u_\epsilon)(\cdot-y_\epsilon)
+Dw_{\delta,k}(u_\epsilon,y_\epsilon,\epsilon)[v,\bar{y}]\ |\ v\in
X_{k},\ \bar{y}
\in \mathbb{R}^{N}\}\nonumber\\
&&+T^{\perp}_{u_\epsilon,y_\epsilon,k}=Y.
\end{eqnarray}
If (\ref{vvcbyryryfg6}) were not true, then there exist
$\epsilon_n\rightarrow 0$ as $n\rightarrow\infty$ such that
$Y_n\neq Y,$ where $Y_n$ denotes the space appeared in the left
side of (\ref{vvcbyryryfg6}) with $\epsilon=\epsilon_n.$ Passing
to a subsequence, we may assume that  $y_{\epsilon_n}\rightarrow
y_k$ and $u_{\epsilon_{_{n}}}\rightarrow u_k$ in $Y$ as
$n\rightarrow\infty,$ since $\{(u_{\epsilon_n},y_{\epsilon_n})\}$
is a bounded sequence in the finite dimensional space $X_k\times
\mathbb{R}^N.$ By  the hypothesis (\ref{hronrenhfbvgyf6rttrf}) and
Proposition \ref{gdgvcfd6d5dttdt}, we deduce that $u_k$ is a
critical point of $I(v+\pi_k(v))$.
  Then by  the conclusion $\bf (iv)$ of Lemma \ref{77eygdjjj876yyyyff5ttt},
$u_k+\pi_k(u_k)$ is a critical point of $I$. We denote it by
$\tilde{u}_k$. Since $D\pi_k(u_k)v\in X$ and
$\mathcal{T}_{u_k}\subset X^\bot$, we get $D\pi_k(u_k)v\bot
\mathcal{T}_{u_k}$, where $\mathcal{T}_{u_k}$ comes from
(\ref{6erwtsfga55a}). Moreover, by Lemma
\ref{77eygdjjj876yyyyff5ttt}, we get that $D\pi_k(u_k)v\in
X^{\bot}_{k}$. Thus,
$$D\pi_k(u_k)v\bot X_k\oplus \mathcal{T}_{u_k}=T_{u_k,0,k}.$$
It follows that the following  subspace of $Y:$
\begin{equation}\label{farae44aea}
\{v-\bar{y}\nabla_x u_k -\bar{y}\nabla_x \pi_k(u_k)+D\pi_k(u_k)v\
|\ v\in X_{k},\ \bar{y} \in \mathbb{R}^{N}\}+ T^\bot_{u_k,0,k}
\end{equation} is equal to
\begin{eqnarray}\label{gs65s5srsra}
&&\{v-\bar{y}\nabla_x u_k -\bar{y}\nabla_x \pi_k(u_k)\ |\ v\in
X_{k},\ \bar{y} \in \mathbb{R}^{N}\}+ T^\bot_{u_k,0,k}\nonumber\\
&=&\{v-\bar{y}\nabla_x \tilde{u}_k \ |\ v\in X_{k},\ \bar{y} \in
\mathbb{R}^{N}\}+ T^\bot_{u_k,0,k}.
\end{eqnarray}
As it has been mentioned above, $\tilde{u}_k=u_k+\pi_k(u_k)\in
\mathcal{K}.$ Therefore, by (\ref{gdtdcfxscsb}), we get that for
every $1\leq j\leq N,$
\begin{eqnarray}\label{hdtd66dtdf}
||\frac{\partial \tilde{u}_k}{\partial
x_j}-\sum^{s}_{i=1}\xi_i(\tilde{u}_k)\frac{\partial u_i}{\partial
x_j}||\leq\sum^{s}_{i=1}\xi_i(\tilde{u}_k) ||\frac{\partial
\tilde{u}_k}{\partial x_j}-\frac{\partial u_i}{\partial
x_j}||\leq\varsigma.
\end{eqnarray}
By $\bf (ii)$ of Lemma \ref{77eygdjjj876yyyyff5ttt} and the fact
that every $\xi_i$ is a Lipschitz function, we deduce  that for
every $1\leq j\leq N,$ as $k\rightarrow\infty$,
\begin{eqnarray}\label{gsgdfrwrwe}
&&||\sum^{s}_{i=1}\xi_i(\tilde{u}_k)\frac{\partial u_i}{\partial
x_j}-\sum^{s}_{i=1}\xi_i(u_k)\frac{\partial u_i}{\partial
x_j}||\nonumber\\
&\leq&\sum^{s}_{i=1}|\xi_i(\tilde{u}_k)-
\xi_i(u_k)|\cdot||\frac{\partial u_i}{\partial x_j}||\leq
C\sum^{s}_{i=1}||\tilde{u}_k- u_k||\cdot||\frac{\partial
u_i}{\partial x_j}||\rightarrow 0,
\end{eqnarray}
where  $C$ is the the Lipschitz constant of $\xi_i.$ By
(\ref{hdtd66dtdf}) and (\ref{gsgdfrwrwe}), we obtain that for
every $1\leq j\leq N,$
$$\limsup_{k\rightarrow\infty}||\frac{\partial \tilde{u}_k}{\partial
x_j}-\sum^{s}_{i=1}\xi_i(u_k)\frac{\partial u_i}{\partial
x_j}||\leq\varsigma.$$ It follows that
\begin{eqnarray}\label{hhdgdg66dtdrer11q}
\limsup_{k\rightarrow\infty}\sup_{|\bar{y}|\leq
1}||\bar{y}\nabla_x
\tilde{u}_k-\sum^{N}_{j=1}\bar{y}_j\sum^{s}_{i=1}\xi_i(u_k)\frac{\partial
u_i}{\partial x_j}||\leq\varsigma.\nonumber
\end{eqnarray}
Thus, when $\varsigma$ is  sufficiently small and $k$
  is  sufficiently large, the
space defined by (\ref{gs65s5srsra}) is equal to $Y$. As a
consequence, when $\varsigma$ is sufficiently small and $k$ is
sufficiently large, the space defined by (\ref{farae44aea}) is
also $Y$. Therefore, the space
\begin{eqnarray}\label{farae4yyy4aea}
&&\{v(\cdot-y_k)-(\bar{y}\nabla_x u_k)(\cdot-y_k)
-(\bar{y}\nabla_x
\pi_k(u_k))(\cdot-y_k)+(D\pi_k(u_k)v)(\cdot-y_k)\nonumber\\
&&\quad \quad |\ v\in X_{k},\ \bar{y} \in \mathbb{R}^{N}\}+
T^\bot_{u_k,y,k}
\end{eqnarray}
is  equal to $Y$. Then we can define a bounded linear operator
\begin{eqnarray}\label{gsts6stsrsf}
&&H_n: Y\rightarrow Y,\nonumber\\
&& w=v(\cdot-y_k)-(\bar{y}\nabla_x u_k)(\cdot-y_k)
-(\bar{y}\nabla_x
\pi_k(u_k))(\cdot-y_k)+(D\pi_k(u_k)v)(\cdot-y_k)+\phi\nonumber\\
&&\mapsto H_n(w)=v(\cdot-y_{\epsilon_{_{n}}})-(\bar{y}\nabla_x
u_{\epsilon_{_{n}}})(\cdot-y_{\epsilon_{_{n}}})
+Dw_{\delta,k}(u_{\epsilon_{_{n}}},y_{\epsilon_{_{n}}},\epsilon_{n})[v,\bar{y}]+\phi,
\nonumber\end{eqnarray} where $\phi\in T^\bot_{u_k,y,k}.$ It
satisfies $Y_n=H_n(Y)$, where $Y_n$ denotes the space appeared in
the left side of (\ref{vvcbyryryfg6}) with $\epsilon=\epsilon_n.$
 By $u_{\epsilon_n}\rightarrow u_k$, $y_{\epsilon_n}\rightarrow
y_k$ and Proposition \ref{gdgvcfd6d5dttdt}, we get that as
$n\rightarrow \infty $,
$$||H_n-id||_{\mathcal{L}(Y)}\rightarrow 0.$$
Therefore,  when  $n$ is large enough, $H_{n}(Y)=Y$. It follows
that $Y_n=Y,$ which  contradicts the assumption. Thus,  when
$k(\delta)$ is large enough and $k\geq k(\delta)$, there exists
$\epsilon_k>0$ such that  if $0\leq\epsilon\leq\epsilon_k,$ then
(\ref{vvcbyryryfg6}) holds. \hfill$\Box$

\section{ Proof of Theorem \ref{bbvv55erdff} }\label{ttrfd444ew3}
By the conclusions $\bf (iii)$ and $\bf (v)$ of Theorem
\ref{77tutyyy},  if $u\in \overline{\mathcal{N}_{\delta,k}}$, then
$\pi_k(u)$ decays exponentially at infinity. Therefore, for $u\in
\overline{\mathcal{N}_{\delta,k}}$ and $y\in \mathbb{R}^N$, we can
define
$$\Gamma_k(u,y)=\int_{\mathbb{R}^{N}}Q_{n^*}(x+y)(u+\pi_k(u))^2dx.$$

By  the same argument as Lemma 3.2 of \cite{Ambrosetti Badiale }
and by (\ref{nvb77dydg8uidjnquzhi}), (\ref{gdbcgtdrrq55drsf}) and
the Lebesgue  Convergence Theorem, we can get the following Lemma:
\begin{lemma}\label{bcvdfg772544} For any given $k\geq k(\delta)$, as $\epsilon\rightarrow
0,$
$$\sup\Big\{
\Big|\frac{1}{\epsilon^{n^*}}\int_{\mathbb{R}^N}V(\epsilon
(x+y))(u+\pi_k(u))^2dx-\Gamma_k(u,y)\Big|\ |\ (u,y)\in
\overline{\mathcal{N}_{\delta,k}}\times\overline{B_{\mathbb{R}^N}(0,R)}\Big\}\rightarrow0$$
and \begin{eqnarray}\label{ggaczddas77ytr} &&\sup\Big\{
\Big|D\Big(\frac{1}{\epsilon^{n^*}}\int_{\mathbb{R}^N}V(\epsilon
(x+y))(u+\pi_k(u))^2dx-\Gamma_k(u,y)\Big)[v,\bar{y}]\Big|\ |\ v\in
X_k,\ ||v||\leq 1, \nonumber\\
&&\quad\quad\quad \bar{y}\in\mathbb{R}^N,\ |\bar{y}|\leq 1,\
(u,y)\in
\overline{\mathcal{N}_{\delta,k}}\times\overline{B_{\mathbb{R}^N}(0,R)}
\Big\}\rightarrow0.\nonumber\end{eqnarray}
\end{lemma}

From now on, for the condition $\bf (V_1)$, we always assume that
 $\triangle Q_{n^*}\geq0$ and $\triangle Q_{n^*}\not\equiv0$  in $
\mathbb{R}^N$, since the proof for the other case is similar.
\begin{lemma}\label{vxbdfreffff} If $\delta>0$ is small enough,
then for any $u\in \overline{\mathcal{N}_{\delta,k}}$,
$\Gamma_k(u,\cdot)$ has a strict local minimum at $y = 0$ and
$D^2_y\Gamma_k(u,0)$ is a positive-definite matrix. More
precisely, there exists a constant $A_k>0$   such that
\begin{eqnarray}\label{bbv88uiiyiyy}
D^2_y\Gamma_k(u,0)y\cdot y\geq A_k|y|^2,\ \forall u\in
\overline{\mathcal{N}_{\delta,k}},\ \forall y\in\mathbb{R}^N.
\end{eqnarray}
\end{lemma}
\noindent{\bf Proof.} By Lemma 4.1 of \cite{Ambrosetti Badiale },
we know that $y=0$ is
 a critical point of $\Gamma_k(u,\cdot)$ for every $u\in \overline{\mathcal{N}_{\delta,k}}$.
 If (\ref{bbv88uiiyiyy}) were not true, then there exist
$\delta_n>0$, $u_n\subset \overline{\mathcal{N}_{\delta_n,k}},$
$n=1,2,\cdots$ and $\{y_n\}\subset S^{N-1}$ such that
$\delta_n\rightarrow 0$ as $n\rightarrow\infty$ and
\begin{eqnarray}\label{vcbfghhfyrt11}
\lim_{n\rightarrow\infty}|D^2_{y}\Gamma_k(u_n,0)y_n\cdot y_n|=0.
\end{eqnarray}
Since $(u_n,y_n)$ is bounded in the  finite dimensional space
$X_k\times \mathbb{R}^N$, passing to a subsequence, we may assume
that $u_n\rightarrow u_0$ in $X_k$,
 and $y_n\rightarrow y_0\in S^{N-1}$ as
$n\rightarrow\infty$. Let $D_{ii}\Gamma_k(u_n,y)$ be the second
derivative of $\Gamma_k(u_n,y)$ with respect to the variable $y_i$
and
$\mbox{diag}\{D_{11}\Gamma_k(u_n,0),\cdots,D_{NN}\Gamma_k(u_n,0)\}$
be diagonal matrix with diagonal elements $D_{11}\Gamma_k(u_n,0),$
$\cdots,$ $D_{NN}\Gamma_k(u_n,0)$. By the appendix of
\cite{Ambrosetti Badiale }, we get that
\begin{eqnarray}\label{yyd990000teook}
D_{ii}\Gamma_k(u_n,0)=-\frac{2}{N}\int_{\mathbb{R}^N}(u_n+\pi_k(u_n))\nabla
Q_{n^*}(x)\cdot  \nabla(u_n+\pi_k(u_n))dx,\ 1\leq i\leq N.
\end{eqnarray}Therefore,
\begin{eqnarray}\label{iiiqttw444343}
D^2_{y}\Gamma_k(u_n,0)y_n\cdot
y_n&=&y^{T}_n\cdot\mbox{diag}\{D_{11}\Gamma_k(u_n,0),\cdots,D_{NN}\Gamma_k(u_n,0)\}\cdot y_n\nonumber\\
&=&-\frac{2}{N}|y_n|^2\int_{\mathbb{R}^N}(u_n+\pi_k(u_n))\nabla
Q_{n^*}(x)\cdot \nabla (u_n+\pi_k(u_n))dx\nonumber\\
&=&-\frac{1}{N}|y_n|^2\int_{\mathbb{R}^N}\nabla Q_{n^*}(x)\cdot
\nabla
(u_n+\pi_k(u_n))^2dx\nonumber\\
&=&\frac{1}{N}|y_n|^2\int_{\mathbb{R}^N}\triangle Q_{n^*}(x)\cdot
(u_n+\pi_k(u_n))^2dx
\end{eqnarray}
By (\ref{vcbfghhfyrt11}) and (\ref{iiiqttw444343}), we infer that
$$\lim_{n\rightarrow\infty}D^2_{y}\Gamma_k(u_n,0)y_n\cdot
y_n=\frac{1}{N}|y_0|^2\int_{\mathbb{R}^N}\triangle Q_{n^*}(x)\cdot
(u_0+\pi_k(u_0))^2dx=0.$$ It is a contradiction, since we have
assumed that $\triangle Q_{n^*}(x)\geq0$  and $\triangle
Q_{n^*}\not\equiv0$ in $ \mathbb{R}^N$. \hfill$\Box$

\medskip

In  the rest of this section, we  assume that $\delta>0$ is
sufficiently small and
  $k\geq k(\delta)$  is sufficiently large such that (\ref{6fgfga4q3q3gvcbvx})
holds, where the constant $k(\delta)$ comes from Theorem
\ref{77tutyyy}.

\bigskip

\noindent{\bf Proof of Theorem \ref{bbvv55erdff}:}

\medskip

By  definition of $\Psi_k(u,y,\epsilon)$ (see
(\ref{hygggvxv555qmmbn})), for $(u,y)\in
\overline{\mathcal{N}_{\delta,k}}\times
\overline{B_{\mathbb{R}^{N}}(0,R)},$
\begin{eqnarray}\label{vxvvdffrrfw44e}
&&\Psi_k(u,y,\epsilon)\nonumber\\
&=&\frac{1}{2}||u(\cdot-y)+w_{\delta,k}(u,y,\epsilon) ||^2+
\frac{1}{2}\int_{\mathbb{R}^N}V(\epsilon x)|u(\cdot-y)+w_{\delta,k}(u,y,\epsilon)|^2dx\nonumber\\
&&-\int_{\mathbb{R}^N}F(u(\cdot-y)+w_{\delta,k}(u,y,\epsilon))dx\nonumber\\
&=&\frac{1}{2}||u(\cdot-y)+w_{\delta,k}(u,y,0) ||^2 +
\frac{1}{2}||w_{\delta,k}(u,y,\epsilon)-w_{\delta,k}(u,y,0)
||^2\nonumber\\
&&+\langle
u(\cdot-y)+w_{\delta,k}(u,y,0),w_{\delta,k}(u,y,\epsilon)-w_{\delta,k}(u,y,0)\rangle\nonumber\\
&&+\frac{1}{2}\int_{\mathbb{R}^N}V(\epsilon
x)|u(\cdot-y)+w_{\delta,k}(u,y,0)|^2dx\nonumber\\
&&+ \frac{1}{2}\int_{\mathbb{R}^N}V(\epsilon
x)|w_{\delta,k}(u,y,\epsilon)-w_{\delta,k}(u,y,0)|^2dx\nonumber\\
&&+\int_{\mathbb{R}^N}V(\epsilon
x)(u(\cdot-y)+w_{\delta,k}(u,y,0))\cdot(w_{\delta,k}(u,y,\epsilon)-w_{\delta,k}(u,y,0))dx\nonumber\\
&&-\int_{\mathbb{R}^N}F(u(\cdot-y)+w_{\delta,k}(u,y,0))dx\nonumber\\
&&-\int_{\mathbb{R}^N}f(u(\cdot-y)+w_{\delta,k}(u,y,0)) \cdot
(w_{\delta,k}(u,y,\epsilon)-w_{\delta,k}(u,y,0))dx\nonumber\\
&&-\eta_1(u,y,\epsilon),
\end{eqnarray}
where
\begin{eqnarray}\label{vvc00kfuufy}
&&\eta_1(u,y,\epsilon)\nonumber\\
&=&\int_{\mathbb{R}^N}F(u(\cdot-y)+w_{\delta,k}(u,y,\epsilon))dx
-\int_{\mathbb{R}^N}F(u(\cdot-y)+w_{\delta,k}(u,y,0))dx\nonumber\\
&&-\int_{\mathbb{R}^N}f(u(\cdot-y)+w_{\delta,k}(u,y,0)) \cdot
(w_{\delta,k}(u,y,\epsilon)-w_{\delta,k}(u,y,0))dx.\nonumber
\end{eqnarray}
By Taylor expansion, we deduce  that there exists
$0<\theta=\theta(x)<1,$ $\forall x\in\mathbb{R}^N$ such that
\begin{eqnarray}\label{jjjs66srrrsnnnzvoozii}
\eta_1(u,y,\epsilon)
&=&\frac{1}{2}\int_{\mathbb{R}^N}f'(u(\cdot-y)+\theta
w_{\delta,k}(u,y,0)+(1-\theta)w_{\delta,k}(u,y,\epsilon))\\
&&\quad\quad\quad
\times(w_{\delta,k}(u,y,\epsilon)-w_{\delta,k}(u,y,0))^2dx\nonumber
\end{eqnarray}
 By the condition $\bf (F_1)$, Proposition \ref{gdgvcfd6d5dttdt} and (\ref{jjjs66srrrsnnnzvoozii}),
 we deduce
that
\begin{eqnarray}\label{fsfdsyyy443wstsr009}
\lim_{\epsilon\rightarrow
0}\sup\{\frac{1}{\epsilon^{n^*}}|\eta_1(u,y,\epsilon) |\ |\
(u,y)\in
\overline{\mathcal{N}_{\delta,k}}\times\overline{B_{\mathbb{R}^N}(0,R)}
\} =0.
\end{eqnarray}
Note that for $v\in X_k,$ $\bar{y}\in\mathbb{R}^N$,
\begin{eqnarray}\label{zhuicundeqiangji}
&&D\eta_1(u,y,\epsilon)[v,\bar{y}]\nonumber\\
&=&\int_{\mathbb{R}^N}f(u(\cdot-y)+w_{\delta,k}(u,y,\epsilon))
\nonumber\\
 &&\quad\quad\times
(v(\cdot-y)-\bar{y}(\nabla_x u)(\cdot-y)+Dw_{\delta,k}(u,y,\epsilon)[v,\bar{y}])dx\nonumber\\
&&-\int_{\mathbb{R}^N}f(u(\cdot-y)+w_{\delta,k}(u,y,0))
\nonumber\\
 &&\quad\quad\times
(v(\cdot-y)-\bar{y}(\nabla_x u)(\cdot-y)+Dw_{\delta,k}(u,y,0)[v,\bar{y}])dx\nonumber\\
&&-\int_{\mathbb{R}^N}f'(u(\cdot-y)+w_{\delta,k}(u,y,0))\cdot(w_{\delta,k}(u,y,\epsilon)-w_{\delta,k}(u,y,0))
\nonumber\\
&&\quad\quad\quad\quad\times (v(\cdot-y)-\bar{y}(\nabla_x
u)(\cdot-y)+Dw_{\delta,k}(u,y,0)[v,\bar{y}])dx\nonumber\\
&&-\int_{\mathbb{R}^N}f(u(\cdot-y)+w_{\delta,k}(u,y,0)) \cdot
(Dw_{\delta,k}(u,y,\epsilon)[v,\bar{y}]-Dw_{\delta,k}(u,y,0)[v,\bar{y}])
\end{eqnarray}
Then by the conclusion $\bf (iii)$ of Theorem \ref{77tutyyy},
Proposition \ref{gdgvcfd6d5dttdt} and the condition $\bf (F_1)$,
we deduce that
\begin{eqnarray}\label{fsfdsyyiiiuy6stsr009}
\lim_{\epsilon\rightarrow
0}\sup\{\frac{1}{\epsilon^{n^*}}||D\eta_1(u,y,\epsilon) ||\ |\
(u,y)\in
\overline{\mathcal{N}_{\delta,k}}\times\overline{B_{\mathbb{R}^N}(0,R)}
\} =0.
\end{eqnarray}
Combining (\ref{fsfdsyyy443wstsr009}) and
(\ref{fsfdsyyiiiuy6stsr009}) yields
\begin{eqnarray}\label{fsfdsyystsr009}
\lim_{\epsilon\rightarrow
0}\sup\{\frac{1}{\epsilon^{n^*}}(|\eta_1(u,y,\epsilon)
|+||D\eta_1(u,y,\epsilon) ||)\ |\ (u,y)\in
\overline{\mathcal{N}_{\delta,k}}\times\overline{B_{\mathbb{R}^N}(0,R)}
\} =0.
\end{eqnarray}

 By the conclusion $\bf (ii)$ of Theorem  \ref{77tutyyy} and the fact that
$$w_{\delta,k}(u,y,\epsilon)-w_{\delta,k}(u,y,0)\in T^{\bot}_{u,y,k},$$
we get
\begin{eqnarray}\label{vcn88tutyty}
&& \langle u(\cdot-y)+w_{\delta,k}(u,y,0)
,w_{\delta,k}(u,y,\epsilon)-w_{\delta,k}(u,y,0)\rangle\nonumber\\
&=&\int_{\mathbb{R}^N}f(u(\cdot-y)+w_{\delta,k}(u,y,0))\cdot
(w_{\delta,k}(u,y,\epsilon)-w_{\delta,k}(u,y,0))dx.
\end{eqnarray}
By Proposition \ref{gdgvcfd6d5dttdt}, we deduce that
\begin{eqnarray}\label{vcbfgt664ta}
&&\eta_2(u,y,\epsilon)\nonumber\\
&:=&\frac{1}{2}||w_{\delta,k}(u,y,\epsilon)-w_{\delta,k}(u,y,0)
||^2+ \frac{1}{2}\int_{\mathbb{R}^N}V(\epsilon
x)|w_{\delta,k}(u,y,\epsilon)-w_{\delta,k}(u,y,0)|^2dx\nonumber\\
&&+\int_{\mathbb{R}^N}V(\epsilon x)(u(\cdot-y)+
w_{\delta,k}(u,y,0))(w_{\delta,k}(u,y,\epsilon)-w_{\delta,k}(u,y,0))dx\nonumber
\end{eqnarray}
also satisfies (\ref{fsfdsyystsr009}). By the conclusion $\bf
(iii)$ of Theorem \ref{77tutyyy}, we infer that
\begin{eqnarray}\label{nvn88tuiyuuy}
J(u(\cdot-y)+w_{\delta,k}(u,y,0))=J(u(\cdot-y)+\pi_k(u)(\cdot-y))=I(u+\pi_k(u)).\end{eqnarray}
Finally, by the conclusions $\bf (iii)$ and $\bf (v)$ of Theorem
\ref{77tutyyy} and (\ref{gdbcgtdrrq55drsf}),  we have
\begin{eqnarray}\label{nv888rtr664}
&&\frac{1}{2}\int_{\mathbb{R}^N}V(\epsilon
x)|u(\cdot-y)+w_{\delta,k}(u,y,0)|^2dx\nonumber\\
&=&\frac{1}{2}\int_{\mathbb{R}^N}V(\epsilon
x)(u(\cdot-y)+\pi_k(u)(\cdot-y))^2dx\nonumber\\
&=&\frac{1}{2}\epsilon^{n^*}\Gamma_k(u,y)+\eta_3(u,y,\epsilon),
\end{eqnarray}
where \begin{eqnarray}\label{nnctf55ftf} \Gamma_k(u,y)&=&
\int_{\mathbb{R}^N}
Q_{n^*}(x)(u(\cdot-y)+\pi_k(u)(\cdot-y))^2dx\nonumber\\
&=&\int_{\mathbb{R}^N} Q_{n^*}(x+y)(u+\pi_{k}(u))^2dx.\nonumber
\end{eqnarray} By
Lemma \ref{bcvdfg772544}, the conclusion $\bf (v)$ of Theorem
\ref{77tutyyy} and (\ref{gdbcgtdrrq55drsf}), we deduce that
$\eta_3$ satisfies (\ref{fsfdsyystsr009}). By
$(\ref{vxvvdffrrfw44e})-(\ref{nv888rtr664})$,  we get that
\begin{eqnarray}\label{vv53rrweeee}
\Psi_k(u,y,\epsilon) =I(u+\pi_k(u))
+\frac{1}{2}\epsilon^{n^*}\Gamma_k(u,y)+\eta(u,y,\epsilon),
\end{eqnarray}
where $\eta=\eta_1+\eta_2+\eta_3$  satisfies
(\ref{fsfdsyystsr009}).

By Lemma \ref{vxbdfreffff},  for every
$u\in\overline{\mathcal{N}_{\delta,k}}$, $\Gamma_k(u,y)$
 has a strict local minimum at $y = 0$ and there is
 a constant $A_k>0$  such that
 \begin{eqnarray}\label{gb8vuvhhhvhvhennab}
 D^2_y \Gamma_k(u,0)\geq A_k \mbox{Id}\end{eqnarray} where Id denotes the
 $N\times N$ identity matrix.
 By
  (\ref{gb8vuvhhhvhvhennab}) and
(\ref{vv53rrweeee}), we deduce that   there exists $\epsilon'_k>0$
such that if $0\leq\epsilon\leq\epsilon'_k$, then
  for every $u\in \overline{\mathcal{N}_{\delta,k}}$, there
exists $y_\epsilon(u)\in B_{\mathbb{R}^N}(0,R/2)$ such that
$y_\epsilon(u)$ is the unique minimizer of
$\Psi_k(u,\cdot,\epsilon)$ in $B_{\mathbb{R}^N}(0,R)$. Moreover,
 by implicit functional theorem,  $y_\epsilon(\cdot)\in
C^1(\overline{\mathcal{N}_{\delta,k}})$. By (\ref{vv53rrweeee}),
we get that
\begin{equation}\label{vvdfrredd33qw}
\lim_{\epsilon\rightarrow
0}||\Psi_k(u,y_\epsilon(u),\epsilon)-I(u+\pi_k(u))
||_{C^1(\overline{\mathcal{N}_{\delta,k}})}=0.
\end{equation}
  By \cite[Theorem IV.3]{CG}, a GM
pair is a special kind of Conley index pair which is associated
with some pseudo-gradient flow of a functional. Therefore, the GM
pair $(W_k,W_k^-)$ which was defined in Remark \ref{uty77ty6bcc}
is a Conley index pair associated with some pseudo-gradient flow
of the functional $g_k(u)=I(u+\pi_k(u)).$ Then by
(\ref{vvdfrredd33qw}) and Theorem III.4 of \cite{CG}, we deduce
that if $\epsilon$ is small enough, then  $(W_k,W_k^-)$ is also a
Conley index pair associated with
 some pseudo-gradient flow of the
functional $\Psi_k(\cdot,y_\epsilon(\cdot),\epsilon).$ By
(\ref{6fgfga4q3q3gvcbvx}) and Theorem 5.5.18 of \cite{cha2}, we
infer that if   $\epsilon$ is sufficiently small, then
  $\Psi_k(\cdot,y_\epsilon(\cdot),\epsilon)$
has at least a critical point $u_{\epsilon}\in
\mathcal{N}_{\delta,k}$. Then by Theorem \ref{ttrgfg77y},
$\tilde{u}_\epsilon:=
u_{\epsilon}(\cdot-y_\epsilon(u_\epsilon))+w_{\delta,k}(u_{\epsilon},y_\epsilon(u_\epsilon),\epsilon)$
is a critical point of $E_\epsilon$. Moreover, by
(\ref{vvdfrredd33qw}), we have
$$\lim_{\epsilon\rightarrow
0}\mbox{dist}_{_{Y}}(\tilde{u}_{\epsilon},\mathcal{K})=0$$ with
$\mathcal{K}=\mathcal{K}^{b}_{a}$. This finishes the proof of
Theorem \ref{bbvv55erdff}. \hfill$\Box$

\section{Appendix A}\label{gdg88futq9oi}
 In this appendix, we shall give the proof of
the existence of $\{\tilde{e}_{j,k}\}$ which satisfies the
conditions $\bf (i)$ and $\bf (ii)$ in Section
\ref{ggdfsrseese11q}.

 Since $X\cap C^{\infty}_{0}(\mathbb{R}^N)$ is dense in $X,$
for any $\mu_k>0,$ we can choose $\{\bar{e}_{j,k}\}\subset X\cap
C^{\infty}_{0}(\mathbb{R}^N)$ such that
\begin{eqnarray}\label{gfh88t7rydgdggg}
\sup_{1\leq j\leq k}||\bar{e}_{j,k}-e'_j||\leq\mu_k\ \mbox{and}\
||\bar{e}_{j,k}||=1,\ 1\leq j\leq k.\end{eqnarray} We show that if
$\mu_k$ is small enough, then $\{\bar{e}_{j,k}\ |\ 1\leq j\leq
k\}\cup\{e_j\ |\ 1\leq j\leq q\}$ is linearly independent. If it
were not true, without loss of generality, we may assume that
\begin{eqnarray}\label{u77d55sfrffdvfd4rr}
\bar{e}_{k,k}=\sum^{k-1}_{j=1}\alpha_j \bar{e}_{j,k}+
\sum^{q}_{j=1}\beta_j e_{j},\end{eqnarray} then
$$\bar{e}_{k,k}=\sum^{k-1}_{j=1}\alpha_j e'_j+
\sum^{k-1}_{j=1}\alpha_j(\bar{e}_{j,k}-e'_j)+
\sum^{q}_{j=1}\beta_j e_{j}.$$ It follows that if
$\mu_k<1/4\sqrt{2}$, then
\begin{eqnarray}\label{gdhbc88fyfg009}
1=||\bar{e}_{k,k}||^2&=&\sum^{k-1}_{j=1}\alpha^2_j+
||\sum^{k-1}_{j=1}\alpha_j(\bar{e}_{j,k}-e'_j)||^2 + 2\langle
\sum^{k-1}_{j=1}\alpha_j e'_j,
\sum^{k-1}_{j=1}\alpha_j(\bar{e}_{j,k}-e'_j)\rangle\nonumber\\
&&+\sum^{q}_{j=1}\beta^2_j + 2\langle \sum^{q}_{j=1}\beta_j e_j,
\sum^{k-1}_{j=1}\alpha_j(\bar{e}_{j,k}-e'_j)\rangle
\nonumber\\
&\geq&
\frac{3}{4}\sum^{k-1}_{j=1}\alpha^2_j+\frac{3}{4}\sum^{q}_{j=1}\beta^2_j
+
||\sum^{k-1}_{j=1}\alpha_j(\bar{e}_{j,k}-e'_j)||^2-8\sum^{k-1}_{j=1}\alpha^2_j||\bar{e}_{j,k}-e'_j||^2
\nonumber\\
&\geq&\frac{1}{2}\sum^{k-1}_{j=1}\alpha^2_j+\frac{1}{2}\sum^{q}_{j=1}\beta^2_j.
\end{eqnarray}
By (\ref{u77d55sfrffdvfd4rr}),
 $$e'_{k}=\sum^{k-1}_{j=1}\alpha_j
e'_j+ \sum^{k-1}_{j=1}\alpha_j(\bar{e}_{j,k}-e'_j)+
\sum^{q}_{j=1}\beta_j e_{j}+(e'_{k}-\bar{e}_{k,k}) ,$$ combining
(\ref{gdhbc88fyfg009}),  we get that
\begin{eqnarray}\label{gdhbc88fdf4yfg009}
1=||e'_k||^2&=&\sum^{k-1}_{j=1}\alpha_j\langle\bar{e}_{j,k}-e'_j,
e'_{k}\rangle+\langle e'_{k}-\bar{e}_{k,k},
e'_{k}\rangle\leq\mu_k\sum^{k-1}_{j=1}|\alpha_j|+\mu_k\nonumber\\
&\leq&(\sqrt{2k}+1)\mu_k.\nonumber
\end{eqnarray}
This induces a contradiction if we assume $(\sqrt{2k}+1)\mu_k<1$.
Thus, $\{\bar{e}_{j,k}\ |\ 1\leq j\leq k\}\cup\{e_j\ |\ 1\leq
j\leq k\}$ is linearly independent if
$\mu_k<\min\{1/(\sqrt{2k}+1),1/4\sqrt{2}\}$.

By (\ref{gfh88t7rydgdggg}) and
$$
\langle \bar{e}_{j,k},\bar{e}_{j',k}\rangle=\langle
e'_j+(\bar{e}_{j,k}-e'_j), e'_{j'}+
(\bar{e}_{j',k}-e'_{j'})\rangle,\ \langle
\bar{e}_{j,k},e_{j'}\rangle=\langle e'_j+(\bar{e}_{j,k}-e'_j),
e_{j'}\rangle,$$ we get that
\begin{eqnarray}\label{gfbb99bjghg}
\sup_{1\leq j,j'\leq k, j\neq j'}|\langle
\bar{e}_{j,k},\bar{e}_{j',k}\rangle|\leq2\mu_k+\mu^2_k,\ \sup_{
j\neq j'}|\langle \bar{e}_{j,k},e_{j'}\rangle|\leq\mu_k.
\end{eqnarray}
Therefore, if $\mu_k$ is   sufficiently small, using Gram-Schmidt
orthogonalizing process to $\{e_j\ |\ 1\leq j\leq
q\}\cup\{\bar{e}_{j,k}\ |\ 1\leq j\leq k\}$,  we get
$\{\tilde{e}_{j,k}\ |\ 1\leq j\leq k\}$ which satisfies the
conditions $\bf (i)$ and $\bf (ii)$ in Section
\ref{ggdfsrseese11q}.

\section{Appendix B}\label{99skdnfhfjfgevfd}
 In this appendix, we give the proof of Proposition
 \ref{gdgvcfd6d5dttdt}.

 Let
\begin{eqnarray}\label{gdgvcb77dydtrr}
\eta_{u,y,k}=(-\triangle+1)^{-1}f(u(\cdot-y)
+\pi_{k}(u)(\cdot-y)).\nonumber
\end{eqnarray}
Then \begin{eqnarray}\label{ggdfvc55rsfrr88u} \eta_{u,y,k}
&=&(-\triangle+1+V(\epsilon x))^{-1}f(u(\cdot-y)
+\pi_{k}(u)(\cdot-y)))\nonumber\\
&&+(-\triangle+1+V(\epsilon x))^{-1}V(\epsilon x)\eta_{u,y,k}.
\end{eqnarray}
Subtracting equation $$S_{u,y,k}\nabla
E_\epsilon(u(\cdot-y)+w_{\delta,k}(u,y,\epsilon))=0$$ from
equation
$$S_{u,y,k}\nabla J(u(\cdot-y)+\pi_{k}(u)(\cdot-y))=0,$$   by
(\ref{ggdfvc55rsfrr88u}) and the mean value theorem, we get that
\begin{eqnarray}\label{fa54qraddassff}
&&L_{u,y,\epsilon,k}\Big(w_{\delta,k}(u,y,\epsilon)-\pi_{k}(u)(\cdot-y)\Big)\nonumber\\
&=&-S_{u,y,k}(-\triangle+1+V(\epsilon x))^{-1}V(\epsilon
x)\eta_{u,y,k}\nonumber\\
&&+S_{u,y,k}(-\triangle+1+V(\epsilon x))^{-1}\Big((f'(
u(\cdot-y)+\tilde{w})-f'(
u(\cdot-y)))\nonumber\\
&&\quad\quad\quad\times(w_{\delta,k}(u,y,\epsilon)-\pi_{k}(u)(\cdot-y))\Big)
\end{eqnarray}
where $\tilde{w}$ lies between $w_{\delta,k}(u,y,\epsilon)$ and
$\pi_{k}(u)(\cdot-y)$. By the conclusion $\bf (iv)$ of Theorem
\ref{77tutyyy}, we get that  $|| w_{\delta,k}(u,y,\epsilon)||\leq
r$ if $0<\delta\leq\delta_r$ and $k\geq k(\delta)$. And by $\bf
(ii)$ of Lemma \ref{77eygdjjj876yyyyff5ttt}, we deduce that if
$k(\delta)$  is large enough and $k\geq k(\delta),$ then
$||\pi_{k}(u)(\cdot-y)||\leq r.$ Therefore, $||\tilde{w}||\leq r$
if  $0<\delta\leq\delta_r$ and $k\geq k(\delta)$. Moreover, by
(\ref{ccsdseee6t6t}), we deduce that if $r$ is small enough,
$0<\delta\leq\delta_r$ and $k\geq k(\delta)$, then
\begin{eqnarray}\label{ggdbvcfdrette0913}
&&\Big|\Big|(-\triangle+1+V(\epsilon x))^{-1}\Big((f'(
u(\cdot-y)+\tilde{w})-f'(
u(\cdot-y))\cdot(w_{\delta,k}(u,y,\epsilon)-\pi_{k}(u)(\cdot-y))\Big)\Big|\Big|\nonumber\\
&&\leq\frac{C}{2}||w_{\delta,k}(u,y,\epsilon)-\pi_{k}(u)(\cdot-y)||,
\end{eqnarray}
where $C$ is the constant appeared in Lemma \ref{yhhdvvvcfd4e}. By
(\ref{ggdbvcfdrette0913}), (\ref{fa54qraddassff}) and Lemma
\ref{yhhdvvvcfd4e},   we get that
\begin{eqnarray}\label{uufgr66etdfvvdf}
C||w_{\delta,k}(u,y,\epsilon)-\pi_{k}(u)(\cdot-y)||\leq
2||(-\triangle+1)^{-1}V(\epsilon x)\eta_{u,y, k}||.
\end{eqnarray}
By (\ref{gdbcgtdrrq55drsf}),  the conclusion $\bf (v)$ of Theorem
\ref{77tutyyy} and \cite[Proposition 4.2]{zhang qi s}, we get that
for any    $n>0$,
\begin{eqnarray}\label{gdbcgtd55drsfooqrq}
&&\sup\{
 ||(1+|x|)^{n}\eta_{u,y, k}||_{L^{\infty}(\mathbb{R}^N)}\ |\
 (u,y)\in\overline{\mathcal{N}_{\delta,k}}\times\overline{B_{\mathbb{R}^N}(0,R)}\}<\infty.
\end{eqnarray}
By (\ref{gdbcgtd55drsfooqrq}), using the same argument as Lemma
3.2 of \cite{Ambrosetti Badiale }, we can get that if $\iota<n^*,$
\begin{eqnarray}\label{cxcv54e3eeew}
\lim_{\epsilon\rightarrow
0}\{\int_{\mathbb{R}^N}\frac{V^2(\epsilon
x)}{\epsilon^{2\iota}}\eta^2_{u,y,k}\ |\
(u,y)\in\overline{\mathcal{N}_{\delta,k}}\times\overline{B_{\mathbb{R}^N}(0,R)}\}=0
\end{eqnarray}
and $$\sup\{\int_{\mathbb{R}^N}\frac{V^2(\epsilon
x)}{\epsilon^{2n^*}}\eta^2_{u,y,k}\ |\
(u,y)\in\overline{\mathcal{N}_{\delta,k}}
\times\overline{B_{\mathbb{R}^N}(0,R)},\
0\leq\epsilon\leq\epsilon^*\}<\infty.$$ Thus, for $\iota<n^*,$
\begin{eqnarray}\label{baobaoainiiii}\lim_{\epsilon\rightarrow
0}\sup\{\frac{1}{\epsilon^\iota}||(-\triangle+1)^{-1}V(\epsilon
x)\eta_{u,y,k}|| \ |\
(u,y)\in\overline{\mathcal{N}_{\delta,k}}\times\overline{B_{\mathbb{R}^N}(0,R)}\}=0
\end{eqnarray}
and
\begin{eqnarray}\label{ggfgttet88etef}
\sup\{\frac{1}{\epsilon^{n^*}}||(-\triangle+1)^{-1}V(\epsilon
x)\eta_{u,y,k}||\ |\
(u,y)\in\overline{\mathcal{N}_{\delta,k}}\times\overline{B_{\mathbb{R}^N}(0,R)},\
0\leq\epsilon\leq\epsilon^*\}<\infty.
\end{eqnarray}
Combining (\ref{uufgr66etdfvvdf}), (\ref{baobaoainiiii}) and
(\ref{ggfgttet88etef}) yields that for $\iota<n^*,$ if $\delta>0$
is small enough and $k\geq k(\delta),$ then
\begin{eqnarray}\label{gdgvcfsf77qrrq}\lim_{\epsilon\rightarrow
0}\{\frac{1}{\epsilon^\iota}||w_{\delta,k}(u,y,\epsilon)-\pi_{k}(u)(\cdot-y)||\
|\
(u,y)\in\overline{\mathcal{N}_{\delta,k}}\times\overline{B_{\mathbb{R}^N}(0,R)}\}=0\end{eqnarray}
and
\begin{eqnarray}\label{ncbcjdgftr745545}
&&\sup\{\frac{1}{\epsilon^{n^*}}||w_{\delta,k}(u,y,\epsilon)-\pi_{k}(u)(\cdot-y)||\
|\
(u,y)\in\overline{\mathcal{N}_{\delta,k}}\times\overline{B_{\mathbb{R}^N}(0,R)},
\ 0\leq\epsilon\leq\epsilon^*\}\nonumber\\
&&<\infty.
\end{eqnarray}

Recall that $S_{u,y,k}: Y\rightarrow T^\bot_{u,y,k}$ is orthogonal
projection. Therefore, for $h\in Y,$
\begin{eqnarray}\label{hhdbcvggdf66tt}
S_{u,y,k}h &=& h-\sum^{q}_{j=1}\langle h, e_j(\cdot-y)\rangle
e_j(\cdot-y)-\sum^{k}_{j=1}\langle h,
\tilde{e}_{j,k}(\cdot-y)\rangle
\tilde{e}_{j,k}(\cdot-y)\nonumber\\
&&-\sum^{N}_{j=1}\Big\langle h,
\sum^{s}_{i=1}\xi_{i}(u)\frac{\partial u_i}{\partial
x_j}(\cdot-y)\Big\rangle
\frac{\sum^{s}_{i=1}\xi_{i}(u)\frac{\partial u_i}{\partial
x_j}(\cdot-y)}{||\sum^{s}_{i=1}\xi_{i}(u)\frac{\partial
u_i}{\partial x_j}||^2}.\nonumber
\end{eqnarray}
Thus, the Fr\'echet partial derivative of $S_{u,y,k}h$  with
respect to $u$ along the vector $v\in X_k$ is
\begin{eqnarray}\label{vvafacczx6trre}
&&D_u(S_{u,y,k}h)[v]\nonumber\\
&=&-\sum^{N}_{j=1}\Big\langle h,
\sum^{s}_{i=1}D\xi_{i}(u)[v]\cdot\frac{\partial u_i}{\partial
x_j}(\cdot-y)\Big\rangle
\frac{\sum^{s}_{i=1}\xi_{i}(u)\frac{\partial u_i}{\partial
x_j}(\cdot-y)}{||\sum^{s}_{i=1}\xi_{i}(u)\frac{\partial
u_i}{\partial x_j}||^2}\nonumber\\
&&-\sum^{N}_{j=1}\Big\langle h,
\sum^{s}_{i=1}\xi_{i}(u)\frac{\partial u_i}{\partial
x_j}(\cdot-y)\Big\rangle
\frac{\sum^{s}_{i=1}(D\xi_{i}(u)[v])\cdot\frac{\partial
u_i}{\partial
x_j}(\cdot-y)}{||\sum^{s}_{i=1}\xi_{i}(u)\frac{\partial
u_i}{\partial x_j}||^2}\nonumber\\
&&+2\sum^{N}_{j=1}\Big(\Big\langle h,
\sum^{s}_{i=1}\xi_{i}(u)\frac{\partial u_i}{\partial
x_j}(\cdot-y)\Big\rangle
\frac{\langle\sum^{s}_{i=1}\xi_{i}(u)\frac{\partial u_i}{\partial
x_j},\sum^{s}_{i=1}(D\xi_{i}(u)[v])\frac{\partial u_i}{\partial
x_j}\rangle}{||\sum^{s}_{i=1}\xi_{i}(u)\frac{\partial
u_i}{\partial x_j}||^4}\nonumber\\
&&\quad\quad\quad\times\sum^{s}_{i=1}\xi_{i}(u)\frac{\partial
u_i}{\partial x_j}(\cdot-y)\Big)
\end{eqnarray}
and the Fr\'echet partial derivative of $S_{u,y,k}h$  with respect
to $y$ along the vector $\bar{y}\in \mathbb{R}^N$ is
\begin{eqnarray}\label{3w33vvafacczx6trre}
&&D_y(S_{u,y,k}h)[\bar{y}]\nonumber\\
&=&\sum^{q}_{j=1}\langle h, (\bar{y}\nabla_x e_j)(\cdot-y)\rangle
e_j(\cdot-y)+ \sum^{k}_{j=1}\langle h,
(\bar{y}\nabla_x\tilde{e}_{j,k})(\cdot-y)\rangle
\tilde{e}_{j,k}(\cdot-y)\nonumber\\
&&+\sum^{q}_{j=1}\langle h, e_j(\cdot-y)\rangle (\bar{y}\nabla_x
e_j)(\cdot-y)+ \sum^{k}_{j=1}\langle h,
\tilde{e}_{j,k}(\cdot-y)\rangle
(\bar{y}\nabla_x\tilde{e}_{j,k})(\cdot-y)\nonumber\\
 &&+\sum^{N}_{j=1}\Big\langle
h, \sum^{s}_{i=1}\xi_{i}(u)\cdot(\bar{y}\nabla_x(\frac{\partial
u_i}{\partial x_j}))(\cdot-y)\Big\rangle
\frac{\sum^{s}_{i=1}\xi_{i}(u)\frac{\partial u_i}{\partial
x_j}(\cdot-y)}{||\sum^{s}_{i=1}\xi_{i}(u)\frac{\partial
u_i}{\partial x_j}||^2}\nonumber\\
&&+\sum^{N}_{j=1}\Big\langle h,
\sum^{s}_{i=1}\xi_{i}(u)\frac{\partial u_i}{\partial
x_j}(\cdot-y)\Big\rangle
\frac{\sum^{s}_{i=1}\xi_{i}(u)\cdot(\bar{y}\nabla_x(\frac{\partial
u_i}{\partial
x_j}))(\cdot-y)}{||\sum^{s}_{i=1}\xi_{i}(u)\frac{\partial
u_i}{\partial x_j}||^2}.
\end{eqnarray}

Differentiating   equations $S_{u,y,k}(\nabla
E_\epsilon(u(\cdot-y)+w_{\delta,k}(u,y,\epsilon)))=0$ and \\
$S_{u,y,k}(\nabla J(u(\cdot-y)+\pi_{k}(u)(\cdot-y))=0$ with
respect to the variable $u$ along the vector $v\in X_k$, we get
that
\begin{eqnarray}\label{11gdgbc66dtfwfwd}
&&S_{u,y,k}(\nabla^2
E_\epsilon(u(\cdot-y)+w_{\delta,k}(u,y,\epsilon))(v(\cdot-y)+Dw_{\delta,k}(u,y,\epsilon)[v,0]))\nonumber\\
&&+D_u(S_{u,y,k}h_1)[v]=0
\end{eqnarray}
and
\begin{eqnarray}\label{gdgbc66dtfwfwd}
&&S_{u,y,k}(\nabla^2
J(u(\cdot-y)+\pi_{k}(u)(\cdot-y))(v(\cdot-y)+D\pi_{k}(u)(\cdot-y)[v,0]))\nonumber\\
&&+D_u(S_{u,y,k}h_2)[v] =0,
\end{eqnarray}
where $h_1=\nabla
E_\epsilon(u(\cdot-y)+w_{\delta,k}(u,y,\epsilon))$ and $h_2=\nabla
J(u(\cdot-y)+\pi_{k}(u)(\cdot-y))$. By (\ref{ggdfvc55rsfrr88u})
and (\ref{ggdbvcfdrette0913}), it is easy to verify that there
exists a constant $C>0$ such that
\begin{eqnarray}\label{yrtegdfdggsfsf111q}
||h_1-h_2||\leq
C||w_{\delta,k}(u,y,\epsilon)-\pi_{k}(u)(\cdot-y)||+C||(-\triangle+1)^{-1}V(\epsilon
x)\eta_{u,y,k}||.
\end{eqnarray}
 By (\ref{yrtegdfdggsfsf111q}) and (\ref{vvafacczx6trre}), we
get that for $||v||\leq 1,$ there exists a constant $C>0$ such
that
\begin{eqnarray}\label{gdttter00oqppp}
&&||D_u(S_{u,y,k}h_2)[v]-D_u(S_{u,y,k}h_1)[v]||\nonumber\\
&\leq&C||w_{\delta,k}(u,y,\epsilon)-\pi_{k}(u)(\cdot-y)||+C||(-\triangle+1)^{-1}V(\epsilon
x)\eta_{u,y,k}||.
\end{eqnarray}
  A direct computation shows that
\begin{eqnarray}\label{bbcb77dtdfff11kkiu}
&&S_{u,y,k}(\nabla^2
E_\epsilon(u(\cdot-y)+w_{\delta,k}(u,y,\epsilon))(v(\cdot-y)+Dw_{\delta,k}(u,y,\epsilon)[v,0]))\nonumber\\
&&-S_{u,y,k}(\nabla^2
J(u(\cdot-y)+\pi_{k}(u)(\cdot-y))(v(\cdot-y)+D\pi_{k}(u)(\cdot-y)[v,0]))\nonumber\\
&=&S_{u,y,k}(\nabla^2
J(u(\cdot-y)+\pi_{k}(u)(\cdot-y))(Dw_{\delta,k}(u,y,\epsilon)[v,0]-D(\pi_{k}(u)(\cdot-y))[v,0]))\nonumber\\
&&-S_{u,y,k}(-\triangle+1)^{-1}\Big\{\Big(
f'(u(\cdot-y)+w_{\delta,k}(u,y,\epsilon))\nonumber\\
&&-
f'(u(\cdot-y)+\pi_{k}(u)(\cdot-y))\Big)\times(v(\cdot-y)+Dw_{\delta,k}(u,y,\epsilon)[v,0])\Big\}
\nonumber\\
&&+S_{u,y,k}(-\triangle+1)^{-1}V(\epsilon
x)\bar{\eta}_{u,y,\epsilon,k}(v)
\end{eqnarray}
where $$ \bar{\eta}_{u,y,\epsilon,k}(v) =(-\triangle+1+V(\epsilon
x))^{-1}(f'(u(\cdot-y)+w_{\delta,k}(u,y,\epsilon)))\cdot(v(\cdot-y)+Dw_{\delta,k}(u,y,\epsilon)[v,0])).
$$
By (\ref{gdbcgtdrrq55drsf}), the conclusion $\bf (v)$ of Theorem
\ref{77tutyyy} and (\ref{l9d8d7f6}) in  $\bf (F_1)$, we get that
for any $v,h\in Y,$ $||v||=||h||=1,$
\begin{eqnarray}\label{qushibadsdsd220}
&&\int_{\mathbb{R}^N}\Big|
f'(u(\cdot-y)+w_{\delta,k}(u,y,\epsilon))-
f'(u(\cdot-y)+\pi_{k}(u)(\cdot-y))\Big|\nonumber\\
&&\quad\quad\quad\times|v(\cdot-y)+Dw_{\delta,k}(u,y,\epsilon)[v,0]|\cdot
|h|dx\nonumber\\
&\leq&
C||w_{\delta,k}(u,y,\epsilon)-\pi_{k}(u)(\cdot-y)||.\nonumber
\end{eqnarray}
It follows that
\begin{eqnarray}\label{gdhf66ey161994}
&&\Big|\Big|(-\triangle+1)^{-1}\Big\{\Big(
f'(u(\cdot-y)+w_{\delta,k}(u,y,\epsilon))\nonumber\\
&&\quad\quad-
f'(u(\cdot-y)+\pi_{k}(u)(\cdot-y))\Big)\times(v(\cdot-y)+Dw_{\delta,k}(u,y,\epsilon)[v,0])\Big\}
\Big|\Big|\nonumber\\
&\leq& C||w_{\delta,k}(u,y,\epsilon)-\pi_{k}(u)(\cdot-y)||.
\end{eqnarray}
By (\ref{11gdgbc66dtfwfwd}), (\ref{gdgbc66dtfwfwd}) and
$(\ref{gdttter00oqppp})-(\ref{gdhf66ey161994})$, we deduce that
\begin{eqnarray}\label{vvcnbd7dyydtdtdt}
&&||S_{u,y,k}(\nabla^2
J(u(\cdot-y)+\pi_{k}(u)(\cdot-y))(Dw_{\delta,k}(u,y,\epsilon)[v,0]-D(\pi_{k}(u)(\cdot-y))[v,0]))||\nonumber\\
&\leq&C||w_{\delta,k}(u,y,\epsilon)-\pi_{k}(u)(\cdot-y)||+C||(-\triangle+1)^{-1}V(\epsilon
x)\eta_{u,y,k}||\nonumber\\
&&+C||(-\triangle+1)^{-1}V(\epsilon
x)\bar{\eta}_{u,y,\epsilon,k}(v)||.\end{eqnarray} By the
conclusion $\bf (ii)$ of Lemma \ref{77eygdjjj876yyyyff5ttt} and
(\ref{ccsdseee6t6t}), we deduce that
\begin{eqnarray}\label{gdgbc77stsrrsf771qa}
&&\lim_{k\rightarrow\infty}\sup\Big\{||\nabla^2
J(u(\cdot-y)+\pi_{k}(u)(\cdot-y))-\nabla^2
J(u(\cdot-y))||_{\mathcal{L}(Y)}\nonumber\\
&&\quad\quad\quad\quad\quad |\
(u,y)\in\overline{\mathcal{N}_{\delta,k}}\times\overline{B_{\mathbb{R}^N}(0,R)}\Big\}=0.\nonumber
\end{eqnarray}
Therefore, as $k\rightarrow\infty$,
\begin{eqnarray}\label{hhfbcg66d5srrsck}
&&||S_{u,y,k}(\nabla^2
J(u(\cdot-y)+\pi_{k}(u)(\cdot-y))(Dw_{\delta,k}(u,y,\epsilon)[v,0]-D(\pi_{k}(u)(\cdot-y))[v,0]))\nonumber\\
&&-S_{u,y,k}(\nabla^2
J(u(\cdot-y))(Dw_{\delta,k}(u,y,\epsilon)[v,0]-D(\pi_{k}(u)(\cdot-y))[v,0]))||\nonumber\\
&&=o(1)||Dw_{\delta,k}(u,y,\epsilon)[v,0]-D(\pi_{k}(u)(\cdot-y))[v,0]||.
\end{eqnarray}
By (\ref{vvcnbd7dyydtdtdt}) and (\ref{hhfbcg66d5srrsck}), we get
that as $k\rightarrow\infty$,
\begin{eqnarray}\label{vvcnbd7duqeeadyydtdtdt}
&&||S_{u,y,k}(\nabla^2
J(u(\cdot-y))(Dw_{\delta,k}(u,y,\epsilon)[v,0]-D(\pi_{k}(u)(\cdot-y))[v,0]))||\nonumber\\
&\leq&C||w_{\delta,k}(u,y,\epsilon)-\pi_{k}(u)(\cdot-y)||+C||(-\triangle+1)^{-1}V(\epsilon
x)\eta_{u,y,k}||\nonumber\\
&&+C||(-\triangle+1)^{-1}V(\epsilon
x)\bar{\eta}_{u,y,\epsilon,k}(v)||\nonumber\\
&&
+o(1)||Dw_{\delta,k}(u,y,\epsilon)[v,0]-D(\pi_{k}(u)(\cdot-y))[v,0]||.\end{eqnarray}
Let $\mathcal{T}_u(\cdot-y)=\{h (\cdot-y)\ |\ h\in
\mathcal{T}_u\}$ and $\mathcal{T}^\bot_u(\cdot-y)$ be the
orthogonal complement space in $Y$,  where $\mathcal{T}_u$ is
defined in (\ref{6erwtsfga55a}). Let
$P_{\mathcal{T}^\bot_u(\cdot-y)}:Y\rightarrow
\mathcal{T}^\bot_u(\cdot-y)$ and
$P_{\mathcal{T}_u(\cdot-y)}:Y\rightarrow \mathcal{T}_u(\cdot-y)$
be orthogonal projections. Since
$Dw_{\delta,k}(u,y,\epsilon)[v,0]\bot X_k(\cdot-y)$ and
$D(\pi_{k}(u)(\cdot-y))[v,0]\bot X_k(\cdot-y),$ where
$X_k(\cdot-y)=\{v(\cdot-y)\ | \ v\in X_k\}$, we deduce that
$$P_{\mathcal{T}^\bot_u(\cdot-y)}(Dw_{\delta,k}(u,y,\epsilon)[v,0]-D(\pi_{k}(u)(\cdot-y))[v,0])\in
T^{\bot}_{u,y,k}.$$
 Therefore, by
Lemma \ref{yhhdvvvcfd4e}, we have
\begin{eqnarray}\label{gbfhvuufydtee}
&&||S_{u,y,k}(\nabla^2
J(u(\cdot-y))P_{\mathcal{T}^\bot_u(\cdot-y)}(Dw_{\delta,k}(u,y,\epsilon)[v,0]
-D(\pi_{k}(u)(\cdot-y))[v,0]))||\nonumber\\
&=&||L_{u,y,0,k}P_{\mathcal{T}^\bot_u(\cdot-y)}(Dw_{\delta,k}(u,y,\epsilon)[v,0]
-D\pi_{k}(u)(\cdot-y))[v,0])||\nonumber\\
 &\geq&
C||P_{\mathcal{T}^\bot_u(\cdot-y)}(Dw_{\delta,k}(u,y,\epsilon)[v,0]-D(\pi_{k}(u)(\cdot-y))[v,0]||.
\end{eqnarray}
Differentiating the following equation  with respect to  variable
$u$ along the vector $v,$
$$\Big\langle w_{\delta,k}(u,y,\epsilon)-\pi_{k}(u)(\cdot-y),
\sum^{s}_{i=1}\xi_i(u)\frac{u_i(\cdot-y)}{\partial
x_j}\Big\rangle=0$$ we get that
\begin{eqnarray}\label{gfb99f0oijfhh}
&&\Big\langle
D(w_{\delta,k}(u,y,\epsilon)-\pi_{k}(u)(\cdot-y))[v,0],
\sum^{s}_{i=1}\xi_i(u)\frac{u_i(\cdot-y)}{\partial
x_j}\Big\rangle\nonumber\\
&=& -\Big\langle w_{\delta,k}(u,y,\epsilon)-\pi_{k}(u)(\cdot-y),
\sum^{s}_{i=1}(D\xi_i(u)[v])\frac{u_i(\cdot-y)}{\partial
x_j}\Big\rangle.\nonumber
\end{eqnarray}
It follows  that there exists a constant $C>0$
 such that\begin{eqnarray}\label{gfbbv8uufy5534}
&&||P_{\mathcal{T}_u(\cdot-y)}(D(w_{\delta,k}(u,y,\epsilon)-\pi_{k}(u)(\cdot-y))[v,0])||\nonumber\\
&\leq& C||w_{\delta,k}(u,y,\epsilon)-\pi_{k}(u)(\cdot-y)||.
\end{eqnarray}
 By
$(\ref{vvcnbd7duqeeadyydtdtdt})-(\ref{gfbbv8uufy5534})$, we deduce
that when $k$ is large enough, then there exists a constant $C>0$
such that
\begin{eqnarray}\label{gfbbv8uudax413qwfy5534}
&&||D(w_{\delta,k}(u,y,\epsilon)-\pi_{k}(u)(\cdot-y))[v,0]||\nonumber\\
&\leq&
C||w_{\delta,k}(u,y,\epsilon)-\pi_{k}(u)(\cdot-y)||+C||(-\triangle+1)^{-1}V(\epsilon
x)\eta_{u,y,\epsilon,k}||.\nonumber\\
&&+C||(-\triangle+1)^{-1}V(\epsilon
x)\bar{\eta}_{u,y,\epsilon,k}(v)||.\nonumber
\end{eqnarray}
Then by $(\ref{baobaoainiiii})- (\ref{ncbcjdgftr745545})$ and the
fact that for $\iota<m,$
\begin{eqnarray}\label{baobaohhainiiii}
&&\lim_{\epsilon\rightarrow
0}\sup\Big\{\frac{1}{\epsilon^\iota}||(-\triangle+1)^{-1}V(\epsilon
x)\bar{\eta}_{u,y,\epsilon,k}(v)||\nonumber\\
&&\quad\quad \quad\quad \ |\
(u,y)\in\overline{\mathcal{N}_{\delta,k}}\times\overline{B_{\mathbb{R}^N}(0,R)},\
v\in X_k,\ ||v||\leq 1\Big\}=0\nonumber
\end{eqnarray}
and
\begin{eqnarray}\label{ggfgtjjtet88etef}
&&\sup\{\frac{1}{\epsilon^{n^*}}||(-\triangle+1)^{-1}V(\epsilon
x)\bar{\eta}_{u,y,\epsilon,k}(v)||\ |\
(u,y)\in\overline{\mathcal{N}_{\delta,k}}\times\overline{B_{\mathbb{R}^N}(0,R)},\nonumber\\
&&\quad\quad\ v\in X_k,\ ||v||\leq 1, \
0\leq\epsilon\leq\epsilon^*\}<\infty,\nonumber\end{eqnarray} we
get that for $\iota<n^*,$
\begin{eqnarray}\label{gdgvcfsf77qrrqhhha}
&&\lim_{\epsilon\rightarrow
0}\sup\Big\{\frac{1}{\epsilon^\iota}||D(w_{\delta,k}(u,y,\epsilon)-\pi_{k}(u)(\cdot-y))[v,0]||\nonumber\\
&&\quad\quad\quad\quad |\
(u,y)\in\overline{\mathcal{N}_{\delta,k}}\times\overline{B_{\mathbb{R}^N}(0,R)},\
v\in X_k,\ ||v||\leq 1\Big\} =0\end{eqnarray} and
\begin{eqnarray}\label{ncbcjdgftr745545hhha}
&&\sup\{\frac{1}{\epsilon^{n^*}}||D(w_{\delta,k}(u,y,\epsilon)-\pi_{k}(u)(\cdot-y))[v,0]||\
|\
(u,y)\in\overline{\mathcal{N}_{\delta,k}}\times\overline{B_{\mathbb{R}^N}(0,R)},\nonumber\\
&&\quad\quad \quad v\in X_k,\ ||v||\leq 1,\
0\leq\epsilon\leq\epsilon^*\}<\infty.
\end{eqnarray}

Differentiating
 the two equations $S_{u,y,k}(\nabla
E_\epsilon(u(\cdot-y)+w_{\delta,k}(u,y,\epsilon)))=0$ and \\
$S_{u,y,k}(\nabla J(u(\cdot-y)+\pi_{k}(u)(\cdot-y))=0$ with
respect to the variable $y$ along the vector $\bar{y}\in
\mathbb{R}^N$, we get that
\begin{eqnarray}\label{gdgbc641wwq6dtfwfwd}
&&S_{u,y,k}(\nabla^2
E_\epsilon(u(\cdot-y)+w_{\delta,k}(u,y,\epsilon))(-\bar{y}\nabla_x
u(\cdot-y)
+Dw_{\delta,k}(u,y,\epsilon)[0,\bar{y}]))\nonumber\\
&&+D_y(S_{u,y,k}h_1)[\bar{y}]=0\nonumber
\end{eqnarray}
and
\begin{eqnarray}\label{gdgbcav66at66dtfwfwd}
&&S_{u,y,k}(\nabla^2
J(u(\cdot-y)+\pi_{k}(u)(\cdot-y))(-\bar{y}\nabla_x u(\cdot-y)
+D(\pi_{k}(u)(\cdot-y))[0,\bar{y}]))\nonumber\\
&&+D_y(S_{u,y,k}h_2)[\bar{y}]=0.\nonumber
\end{eqnarray}
The same arguments as (\ref{gdgvcfsf77qrrqhhha}) and
(\ref{ncbcjdgftr745545hhha}) yield that for $\iota<n^*,$
\begin{eqnarray}\label{qwasgdgvcfsf77qrrqhhha}
&&\lim_{\epsilon\rightarrow
0}\sup\Big\{\frac{1}{\epsilon^\iota}||D(w_{\delta,k}(u,y,\epsilon)-\pi_{k}(u)(\cdot-y))[0,\bar{y}]||\nonumber\\
&&\quad\quad\quad\quad\ |\
(u,y)\in\overline{\mathcal{N}_{\delta,k}}\times\overline{B_{\mathbb{R}^N}(0,R)},\
\bar{y}\in \mathbb{R}^N,\ |\bar{y}|\leq 1\Big\}
=0\nonumber\end{eqnarray}
 and
\begin{eqnarray}\label{qeasncbcjdgftr745545hhha}
&&\sup\{\frac{1}{\epsilon^{n^*}}||D(w_{\delta,k}(u,y,\epsilon)-\pi_{k}(u)(\cdot-y))[0,\bar{y}]||\
|\
(u,y)\in\overline{\mathcal{N}_{\delta,k}}\times\overline{B_{\mathbb{R}^N}(0,R)},\nonumber\\
&&\quad\quad \quad \bar{y}\in \mathbb{R}^N,\ |\bar{y}|\leq 1,\
0\leq\epsilon\leq\epsilon^*\}<\infty.\nonumber
\end{eqnarray}

\bigskip

\noindent{\bf Acknowledgements} The author would like to thank the
referee  for her or his comments and suggestions on the
manuscript. This work was supported by NSFC (10901112) and BNSF
(1102013).


\begin{thebibliography}{99}

\bibitem{Ambrosetti Badiale } A. Ambrosetti, M. Badiale,
S. Cingolani,  {\it Semiclassical states of nonlinear
Schr\"odinger equations}, Arch. Ration. Mech. Anal.  {\bf 140}
(1997) 285-300.

\bibitem{Ambrosetti Secchi} A. Ambrosetti, A. Malchiodi, S.
Secchi,  {\it Multiplicity results for some nonlinear
Schr\"odinger equations with potentials}, Arch. Ration. Mech.
Anal. {\bf  159} (2001) 253-271.

\bibitem{Ambrosetti peral}
A. Ambrosetti, J. Garcia Azorero, I. Peral,  {\it Perturbation of
$\triangle u+u^{(N+2)/(N-2)}=0$, the Scalar Curvature Problem in
$\mathbb{R}^N$, and Related Topics}, J. Funct. Anal. {\bf  165}
(1999) 117-149.

\bibitem{Ber-Lions}
H. Berestycki, P.L. Lions,  {\it Nonlinear scalar field equations:
I, II}, Arch. Ration. Mech. Anal. {\bf 82} (1983) 313-375.


\bibitem{Beyon-Jeanjean}J. Byeon, L. Jeanjean,
{\it Standing waves for nonlinear Schr\"oinger equations with a
general nonlinearity}, Arch. Ration. Mech. Anal. {\bf 185} (2007)
185-200.

\bibitem{Beyon-Jeanjean2}J. Byeon, L. Jeanjean,  K. Tanaka,
{\it Standing waves for nonlinear Schr\"oinger equations with a
general nonlinearity: One and two dimensional cases}, Comm.
Partial Differential Equations {\bf 33} (2008), 1113-1136.


\bibitem{cha}K.- C. Chang, {\it Infinite Dimentional Morse Theory and
Multiple Solution Problem}, Birkh\"auser Boston Inc., Boston, MA,
1991.

\bibitem{cha2}K.- C. Chang, {\it  Methods in Nonlinear Analysis,}
Springer-Verlag, Berlin, Heidelberg, 2005.

\bibitem{CG}K.- C. Chang, N. Ghoussoub,
{\it The Conley index and the critical groups via an extension of
Gromoll-Meyer theory}, Topol. Methods Nonlinear Anal. {\bf 7}
(1996)  77-93.

\bibitem{jeanjean3}S. Cingolani, L. Jeanjean, S. Secchi,
{\it Multi-peak solutions for magnetic NLS equations without
non-degeneracy conditions}, ESAIM Control Optim. Calc. Var. {\bf
15} (2009) 653-673.


\bibitem{Dancer} E. N. Dancer,  {\it On the uniqueness of the positive
solution of a singularly perturbed problem}, Rocky Mountain J.
Math. {\bf 25} (1995) 957-975.



\bibitem{pino} M. Del Pino, P. Felmer,
{\it Local mountain passes for semilinear elliptic problems in
unbounded domains}, Calc. Var. Partial Differential Equations {\bf
4} (1996) 121-137.

\bibitem{Del felmer} M. Del Pino, P. Felmer,
{\it Semi-classical states of nonlinear Schr\"odinger equations: a
variational reduction method}, Math. Ann. {\bf 324} (2002) 1-32.

\bibitem{Del felmer2}M. Del Pino,  P. Felmer, {\it Semi-classical states for nonlinear
Schr\"odinger equations}, J. Funct. Anal. {\bf 149} (1997)
245-265.

\bibitem{Del felmer3}
M. Del Pino,  P. Felmer, {\it Multi-peak bound states for
nonlinear Schr\"odinger equations}, Ann. Inst. H. Poincar\'e Anal.
Non Lin\'eaire  {\bf 15} (1998)  127-149.

\bibitem{Grossi} M. Grossi,
{\it Some results on a class of nonlinear Schr\"odinger
equations}, Math. Z.  {\bf 235} (2000)  687-705.


\bibitem{Gui} C. Gui,
{\it Existence of multi-bump solutions for nonlinear Schr\"odinger
equations}, Comm. Partial Differential Equations  {\bf 21} (1996)
787-820.




\bibitem{Felli}
V. Felli, M. Schneider,  {\it Perturbation results of critical
elliptic equations of Caffarelli-Kohn-Nirenberg type},  J.
Differential Equations  {\bf 191} (2003)  121-142.


\bibitem{GM}D. Gromoll, W. Meyer,  {\it On differentiable functions with isolated critical
points},  Topology  {\bf 8} (1969)  361-369.

\bibitem{KW}X. Kang, J. C. Wei,
{\it On interacting bumps of semi-classical states of nonlinear
Schr\"odinger equations},  Adv. Differential Equations  {\bf 5}
(2000)  899-928.

 \bibitem{Rabi2}
P. H. Rabinowitz, {\it On a class of nonlinear Schr\"odinger
equations}, Z. Angew. Math. Phys.  {\bf 43} (1992)  270-291.

\bibitem{Simon} B. Simon,  {\it Schr\"odinger semigroups},
 Bull. Amer. Math. Soc. (N.S.)  {\bf 7} (1982)  447-526.

\bibitem{spanier} E. H. Spanier, {\it  Algebraic Topology},
Springer-Verlag, New York,   1966.





\bibitem{Will}M. Willem, {\it  Minimax Theorems},
Progress in Nonlinear Differential Equations and their
Applications  24. Birkh\"auser Boston, Inc., Boston, MA, 1996.

\bibitem{zhang qi s} Qi S. Zhang, {\it Positive solutions to  $\triangle
u-Vu+Wu^p=0$ and its parabolic counterpart in noncompact
manifolds},  Pacific J. Math.  {\bf 213} (2004)  163-200.


\end{thebibliography}
\end{document}